\documentclass[final,3p,times]{elsarticle}
\usepackage{lineno,hyperref}
\usepackage{amssymb}
\usepackage{amsmath}
\usepackage{float}
\usepackage{epstopdf}
\usepackage{subfigure}
\usepackage{amsthm}
\usepackage{ifthen}
\numberwithin{equation}{section}
\modulolinenumbers[1]
\usepackage{graphics}
\usepackage{graphicx} 
\newcommand{\matlab}{\text{\includegraphics[height=0.8em]{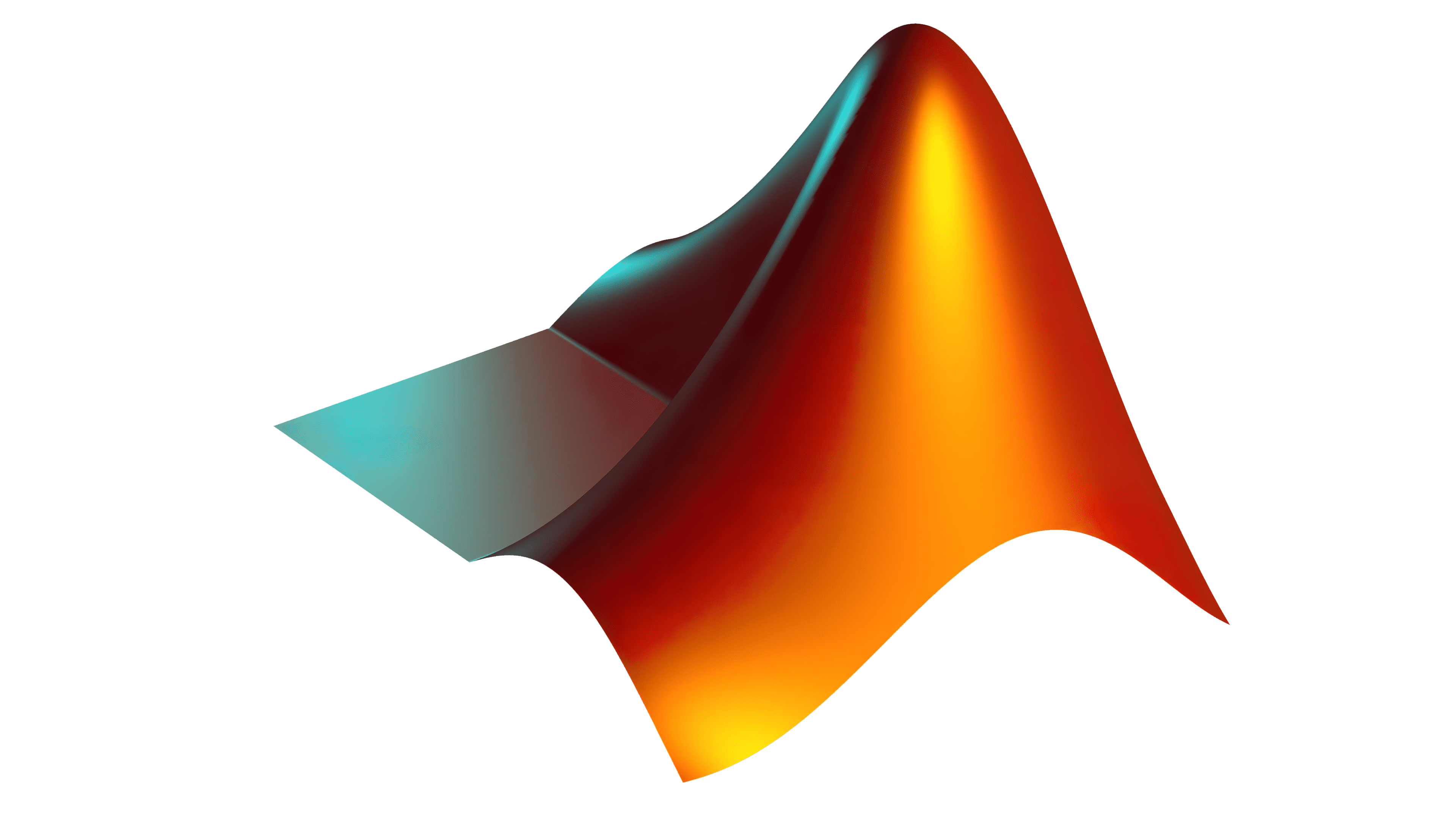}}}
\usepackage[mathscr]{eucal}
\usepackage{hyperref}
\usepackage[utf8]{inputenc}
\usepackage{array,multirow}
\usepackage{comment}
\usepackage{mathtools}
\usepackage{eucal}
\usepackage{nicematrix}
\usepackage{times}
\usepackage[english]{babel}

\usepackage{fourier} 
\usepackage{array}
\usepackage{makecell}
\usepackage{tikz}

\newcommand{\circnum}[1]{\tikz[baseline=(char.base)]\node[draw,circle,inner sep=1pt](char){#1};}

\usepackage[utf8]{inputenc}
\usepackage[english]{babel}
\usepackage{booktabs} 
\usepackage{graphicx} 
\usepackage{enumitem} 
\usepackage{tikz}
\usepackage{multicol}
\usepackage{float}
\usepackage{color} 
\definecolor{mygreen}{RGB}{28,172,0} 
\definecolor{mylilas}{RGB}{170,55,241}

\usepackage{graphicx}
\usepackage{pstricks}
\usepackage{comment}
\usepackage{listings}
\usepackage{xcolor}
\definecolor{codegreen}{rgb}{0,0.6,0}
\definecolor{codegray}{rgb}{0.5,0.5,0.5}
\definecolor{codepurple}{rgb}{0.58,0,0.82}
\definecolor{backcolour}{rgb}{0.95,0.95,0.92}
\definecolor{corn}{rgb}{0.98, 0.93, 0.36}
\definecolor{pastelyellow}{rgb}{0.99, 0.99, 0.9}
\definecolor{red-violet}{rgb}{0.78, 0.08, 0.52}
\lstdefinestyle{matlab}{
    backgroundcolor=\color{backcolour},
    commentstyle=\color{codegreen},
    keywordstyle=\color{magenta},
    numberstyle=\tiny\color{black},
    stringstyle=\color{codepurple},
    basicstyle=\ttfamily\footnotesize,
    breakatwhitespace=false,
    breaklines=true,
    captionpos=b,
    keepspaces=true,
    numbers=left,
    numbersep=5pt,
    showspaces=false,
    showstringspaces=false,
    showtabs=false,
    tabsize=2
}

\usepackage{matlab-prettifier}
\usepackage{algorithm}
\usepackage{algpseudocode}
\usepackage{algorithmicx}
\usepackage{caption}
\usepackage{subcaption}

\newcommand{\md}{{\rm d}}
\newcommand{\dx}{\, \md \textbf{x}}

\newcommand{\dgamma}{\, \md \gamma}
\newcommand{\grad}{\nabla}

\newtheorem{thm}{Theorem}[section]

\theoremstyle{definition}

\newtheorem{example}[thm]{Example}
\usepackage{multicol}
\usepackage{color} 
\usepackage{etoolbox}
\makeatletter
\renewcommand{\ps@pprintTitle}{}
\makeatother

\newcounter{remarkcounter}
\newcommand{\remark}[1]{%
    \refstepcounter{remarkcounter}
    \vskip .3cm%
    \noindent%
    \textbf{Remark \theremarkcounter.}~\textit{#1}\par%
}

\newcommand{\note}[1]{\vskip .3cm \noindent \textbf{Note:}~\textit{#1}\par}

\begin{document}
\begin{frontmatter}



\title{Vectorized 3D mesh refinement and implementation \\ of primal hybrid FEM in MATLAB}


\author[1]{Harish Nagula Mallesham 
}\ead{mat21h.nagula@stuiocb.ictmumbai.edu.in}

\author[2]{Sharat Gaddam 
}\ead{gaddamsharat@gmail.com}

\author[3,4]{Jan Valdman 
}\ead{ jan.valdman@utia.cas.cz }

\author[1]{Sanjib Kumar Acharya \corref{cor1}
}\ead{sk.acharya@iocb.ictmumbai.edu.in}

 \affiliation[1]{organization={Institute of Chemical Technology Mumbai IndianOil Odisha Campus},
 city={Bhubaneswar, Odisha}, postcode={751013},
 country={India}}
 
 \affiliation[2]{organization={High Energy Materials Research Laboratory},
 city={Sutarwadi, Pune, Maharashtra}, postcode={411021},
 country={India}}

\affiliation[3]{organization={Faculty of Information Technology, Czech Technical University in Prague},
 addressline={Thákurova 9, 16000 Prague},
 country={Czech Republic}}

 \affiliation[4]{organization={The Czech Academy of Sciences, Institute of Information Theory and Automation},
 city={Pod vod\'{a}renskou v\v{e}\v{z}\'{i} 4, 18208, Prague 8},
 country={Czech Republic}}
\cortext[cor1]{Corresponding Author}

\begin{abstract}
     
    In this article, we introduce a \textit{Face-to-Tetrahedron} connectivity in MATLAB together with a vectorized 3D uniform mesh refinement technique. We introduce a MATLAB vectorized assembly of 3D lowest-order primal hybrid finite element matrices for a second-order elliptic problem. We introduce a parallel solver and a vectorized Schur complement solver to solve the associated linear problem. 
     The numerical results illustrate the software's runtime performance.
\end{abstract}

\begin{keyword}
   3D mesh refinement, finite elements, primal hybrid, Schur complement, MATLAB, vectorization


\end{keyword}

\end{frontmatter}

\section{Introduction}

Using "for loops" to assemble finite element matrices and refine meshes are the two most time-consuming MATLAB operations. 
In MATLAB, there aren't many publications that cover quick mesh refinement and assembly methods. 
A rapid construction process for the stiffness matrix was initially presented by Rahman and Valdman \cite{rahman2013fast}, who vectorized the for loops utilizing matrix-wise array operations in 2D and 3D. The vectorized assembly of several FE matrices is then presented in a few articles, including \cite{anjam2015fast,moskovka2022fast,cuvelier2016efficient,funken2011efficient,moskovka2024vectorizedbasiclinearalgebra}. 
 In arbitrary spatial dimensions, Feifel and Funken \cite{feifel2024efficient} offer a vectorized assembly of stiffness matrix and adaptive mesh refinement based on the newest vertex bisection. The authors of \cite{bey2000simplicial} talk about 3D red-refinement using Freudenthal's technique, which splits each tetrahedron into eight smaller tetrahedra.
 
  We present a vectorized MATLAB 3D uniform mesh refinement using \textit{ Nodes-to-Edge} connectivity and subdividing each tetrahedron into 12 subtetrahedra.
 In implementing 3D FEM where the degrees of freedom of finite elements are face oriented, it is important to know the connectivity of faces to the tetrahedron for computing the face integrals. 
 We use \textit{Nodes-to-Edge} and \textit{Edges-to-Tetrahedron} connectivity, which will establish \textit{Face-to-Tetrahedron} connectivity. We remark here that instead of using the nodes of the face to develop \textit{Face-to-Tetrahedron} connectivity \cite{feifel2024efficient}, we use the edges of the face to generate \textit{Face-to-Tetrahedron} connectivity.

In \cite{mallesham2025vectorized}, the authors present a vectorized implementation of the lowest-order 2D primal hybrid FEM in MATLAB. This FEM was originally introduced by Raviart and Thomas in 1977 (\cite{raviart}) for second-order elliptic problems, for which optimal-order \textit{a priori} error estimates have been established. 
This method allows for the simultaneous computation of two unknown variables (primal and hybrid) without losing accuracy, which is an advantage in many physical situations.Its finite element functions are discontinuous across the interelement boundaries in the mesh, and the weak continuity is imposed by the Lagrange multiplier. Thus, we can compute the primal variable in parallel by solving the Schur complement system for the Lagrange multiplier. For a review of the literature, see \cite{Quarteroni1979,milner,park,flux,Bendali2025,acharya, kamna, Aposteriori, acharyaElasticity,ACHARYA202516}.       

 We present a vectorized MATLAB implementation of the lowest-order 3D primal hybrid FEM for second-order elliptic problems with mixed boundary conditions.

 Our major contributions are as follows.

\begin{itemize}
\item A vectorized 3D red-refinement strategy and a vectorized \textit{Nodes-to-Edge}, \textit{Edges-to-Tetrahedron}, and \textit{Face-to-Tetrahedron} connectivity.
\item  Vectorized assembly of the lowest-order 3D primal hybrid FE matrices for the second-order problem.
\item  A parallel and a vectorized Schur complement solver, which are significantly faster than the MATLAB linear solver.
\item The software is available for download through the link: \newline
 \matlab~\url{https://in.mathworks.com/matlabcentral/fileexchange/181928}
\end{itemize}
 Our approach to modeling the problem is shorter (see \eqref{eq1}) and can be easily adapted to modified problems, such as those that involve convection and non-linear terms.
The numerical experiments were carried out using MATLAB R2025a on a computer equipped with a $\mathrm{i}5-10210\mathrm{U}$ CPU operating at a frequency of 1.60GHz-2.10 GHz, with 16GB  RAM and a $\mathrm{x}64$-based processor with 1TB of system memory.

In Section \ref{sub:Data structureRefinement}, we define the orientation of tetrahedrons and their faces. In addition, we define the initial data structure for implementations. In Section \ref{sub:connectivity}, we build the connectivity of the \textit{Nodes-to-Edge}, \textit{Edges-to-Tetrahedron}, and \textit{Faces-to-Tetrahedron}; after that, details of the red-refinement on a tetrahedron. Section \ref{sub:primalhybrid} contains the primal hybrid FEM formulation for the elliptic problem and the MATLAB assembly of the stiffness, mass, and the Lagrange multiplier matrix. Additionally, it includes the implementation of a parallel and a vectorized Schur complement algorithm. We give our concluding remarks in Section \ref{sub:conclusion}.

\section{Data structure and orientation}\label{sub:Data structureRefinement}
Let $\Omega=(0,1)^3$ be the computational domain with boundary $\Gamma=\overline{\Gamma}_{D}\cup\overline{\Gamma}_N$, where $\overline{\Gamma}_D$  and $\Gamma_N=\Gamma/\overline{\Gamma}_{D}$ are the Dirichlet and Neumann boundaries, respectively.  Let $\mathcal{T}_{h}$ be a regular family of triangulations of the set $\overline{\Omega}$  with shape regular tetrahedrons $T$ whose diameters are $\leq h$ such that
\begin{equation}
\overline{\Omega}=\bigcup_{T\in\mathcal{T}_{h}}\overline{T}.\nonumber
\end{equation}

\textbf{Notations:} 

\noindent\begin{minipage}{.5\linewidth}
\begin{align*}
\mathcal{F}_h &= \text{set of all 2-simplicies in} ~\mathcal{T}_h,
  \\
\mathcal{F}^h_{\Omega} &= \text{ set of all interior faces of}~\mathcal{F}_h, \\
\mathcal{F}^h_{D} &= \text{set of all Dirichlet boundary faces of}~\mathcal{F}_h, \\
\mathcal{F}^h_{N} &= \text{set of all Neumann boundary faces of}~\mathcal{F}_h,\\
\mathcal{E}_h &= \text{set of all edges in}~\mathcal{T}_h \\
\mathcal{V}_h & = \text{set of all nodes/vertices in}~\mathcal{T}_h, \\
Ct&=\text{centroid of face},\\
\end{align*}
\end{minipage}%
\begin{minipage}{.5\linewidth}
\begin{align*}
\overrightarrow{pq}= \text{vector}~&\text{joining initial point}~p~\text{and final point}~ q,\\
\overrightarrow{ab} \times \overrightarrow{cd}= \text{is} ~& \text{the cross product of vectors}~\overrightarrow{ab}~\text{and}~\overrightarrow{cd},\\
n(A)=\text{num}&\text{ber of elements in A (vector or matrix)},\\
\texttt{nE}= n(\mathcal{T}_h),&\\
\texttt{nC}=n(\mathcal{V}_h),~& \\
\texttt{nF}=n(\mathcal{F}_h),& \\
\texttt{nEd}=n(\mathcal{E}_h)&. \\
\end{align*}
\end{minipage}


\begin{figure}[H]             
\centering
\includegraphics[height=7cm]{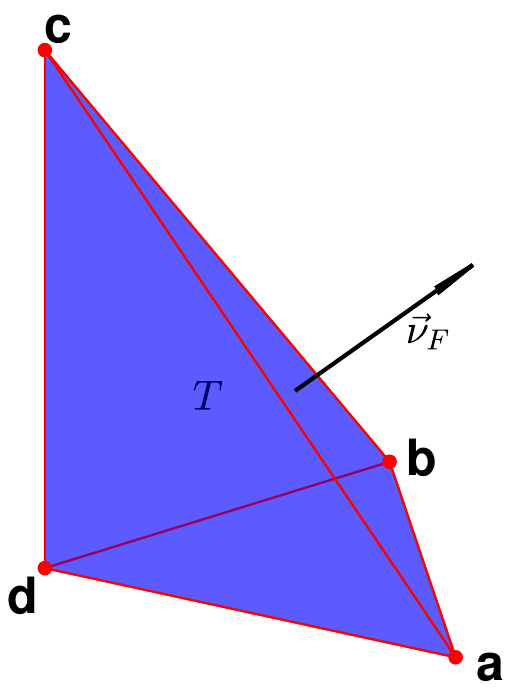}
\caption{\footnotesize  Tetrahedron $T=[a~~ b ~~c~~d]$ with an outward unit normal $\vec{\nu}_F$ to a face $F$.}
 \label{fig:3DTetrahedron}
 \end{figure}

\subsection{Face orientation}\label{sub:FaceOrientation}
The 3-tuple represented by $[x~y~z]$ is an oriented face $F$ of the tetrahedron $T$, denoted by the 4-tuple $[a~b~c~d]$ as in \figurename{~\ref{fig:3DTetrahedron}}, if $\overrightarrow{yz} \times  \overrightarrow{yx}$ is an outward normal to the tetrahedron $T$ for $x,y,z \in \{a,b,c,d\}$.

\noindent\textbf{Note}: We say that the two oriented faces $F_1=[x_1~y_1~z_1]$ and $F_2=[x_2~y_2~z_2]$ of the tetrahedron $T=[a~b~c~d]$ are the same if $\{x_1,y_1,z_1\}=\{x_2,y_2,z_2\}$. For example, \figurename{~\ref{fig:3DTetrahedron}} can be reproduced by the script $$\texttt{fig\_TetrahedronOutwardUnitNormal.m}$$
which shows face $F=[a~b~c]=[b~c~a]=[c~a~b]$.  

\subsection{Tetrahedron orientation}
\begin{figure}[H]             
\centering
\includegraphics[height=7cm]{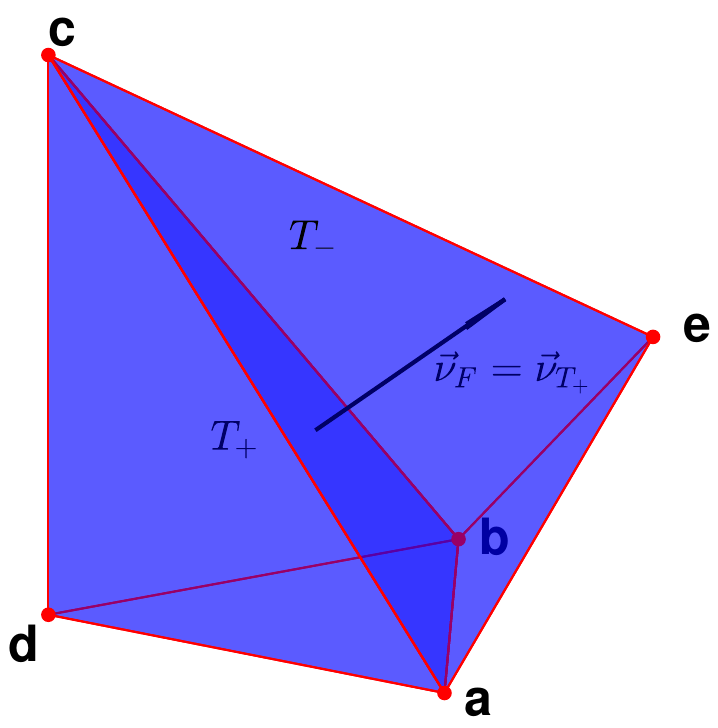}
\caption{Tetrahedrons $T_{+}$ and $T_{-}$ sharing face $[a~~ b ~~c]$.}
 \label{fig:3DTetrahedronsOrientation}
 \end{figure}

As defined in Subsection \ref{sub:FaceOrientation}, a 4-tuple $T=[a~b~c~d]$ always represents a tetrahedron. This representation should always make sure that the first 3-tuple
of $T$, that is, $[a ~b ~c]$, is the oriented face $F$ of the tetrahedron $T$.
List of all oriented faces of the tetrahedron $[a~ b ~c ~d]$ are $[a ~b~ c],~[a~ c ~d],~[a ~d ~b],$ and $[d~ c~ b]$. In this case, we say that the tetrahedron $T$ is in orientation with all four faces.
For each of these four faces, there exists a tetrahedron in orientation and a tetrahedron in opposite orientation( only in the case of an interior face). The \figurename{~\ref{fig:3DTetrahedronsOrientation}} can be reproduced by the script
$$\texttt{fig\_TetrahedronsOrientation.m}$$
which depicts the tetrahedron in the opposite orientation to the face $F=[a~b~c]$, where $T_{-}$ is in the opposite orientation to the face $F$.



 Let $F\in\mathcal{F}^{h}_{\Omega}$ be  such that
$F=\partial T_{+}\cap\partial T_{-}$  where $T_{+}$ and $T_{-}$ are adjacent triangles in $\mathcal{T}_{h}$ (as shown in \figurename{~\ref{fig:3DTetrahedronsOrientation}}) with $\vec{\nu}_{T_+}=\vec{\nu}_F$ denoting the unit normal of $F$ pointing from $T_+$ to $T_-$ and $\vec{\nu}_{T_-}=-\vec{\nu}_F$. We define the jump of a scalar-valued function $v$ on the face $F$ as
\begin{align*}
{[\![v]\!]}_{F}={({v|}_{T_{+}})|}_{F}-{({v|}_{T_{-}})|}_{F}.
\end{align*} 
Further, for $F\in\mathcal{F}^h_D$ such that
 $F=\partial T\cap\overline{\Gamma}_D$ for some $T\in\mathcal{T}_h$, we define ${[\![v]\!]}_{F}=({v|}_{T})|_{F}$. \newline
\textbf{Note}: A face is simply a 2-simplex shared by two tetrahedrons or the boundary of a tetrahedron.
Each face has two distinct oriented faces with respect to the outward normal directions. For example, look at \figurename{~\ref{fig:3DTetrahedron}}; the face with nodes $\{ a,b,c\}$ has two oriented faces $[a~b~c]$ and $[c~b~a]$, by using Subsection \ref{sub:FaceOrientation} we can say that $[a~b~c]$ is in orientation with $T=[a~b~c~d]$ whereas the latter $[c~b~a]$ is not.


\lstset{language=Matlab,numbers=left,
style= Matlab-editor,
backgroundcolor=\color{pastelyellow},
basicstyle=\ttfamily\footnotesize,
commentstyle=\color{codegreen},
numberstyle=\tiny\color{black}, 
stringstyle=\color{red-violet},
breakatwhitespace=false,
breaklines=true,
    captionpos=b,
    keepspaces=true,
    numbers=left,
    numbersep=5pt,
    showspaces=false,
    showstringspaces=false,
    showtabs=false,
    tabsize=2}

\subsection{Initial data structure}
\begin{figure}[H]             
\centering
\vspace{-0.45cm}
\hspace{-1.8cm}
\includegraphics[height=7cm]{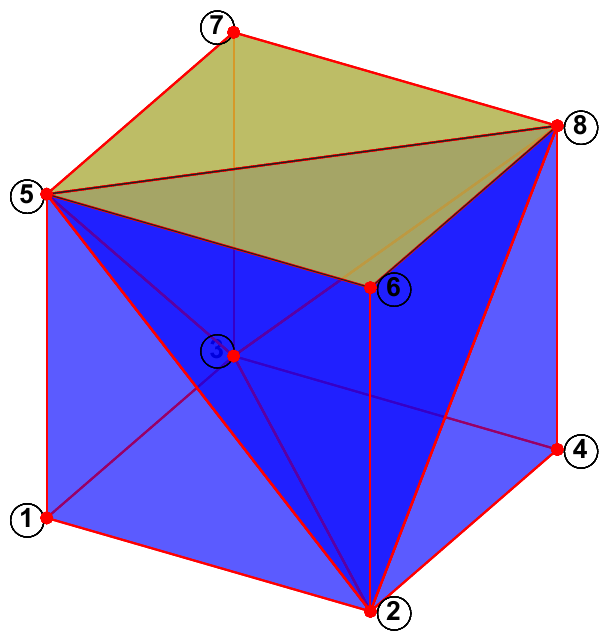}
\caption{\footnotesize  Tetrahedral mesh of a unit cube composed of five tetrahedron. The two yellow boundary faces represent $\mathcal{F}^h_{D}$, and the remaining boundary faces represent $\mathcal{F}^h_{N}$.
}
 \label{fig:3DInitialMesh}
 \end{figure}

For the initial triangulation, we consider the cube $\overline{\Omega}$. The coordinates, elements, Dirichlet boundary, and Neumann boundary information for the initial mesh are represented by the matrices \texttt{Coord}, \texttt{Tetra}, \texttt{Db}, and \texttt{Nb}, respectively, in Listing \ref{lst:InitMeshData}.
\figurename{~\ref{fig:3DInitialMesh}} shows an example of the initial mesh with five tetrahedra, generated from the matrices below. 
\begin{lstlisting}[numbers=none,caption= Initial mesh data, label=lst:InitMeshData]
Coord=[0 0 0; 1 0 0; 0 1 0; 1 1 0; 0 0 1; 1 0 1; 0 1 1; 1 1 1]; 
Tetra=[1 2 3 5; 5 8 6 2; 7 3 8 5; 4 8 3 2; 8 3 2 5];
Db=[7 8 5; 5 8 6];
Nb=[1 5 2; 5 6 2; 1 3 5; 7 5 3; 7 3 8; 4 8 3;  2 6 8; 4 2 8; 4 3 2; 1 2 3];
\end{lstlisting}


\section{Connectivity and refinement}\label{sub:connectivity}
In this section, first we discuss the connectivity of \textit{Nodes-Edges-Tetrahedrons} through functions \textbf{NumEdges} and \textbf{EdgeNum2Tetra}.
\subsection{Nodes-to-Edge Connectivity} 
The function \textbf{NumEdges} establishes that 
\textit{Node-to-Edge} connectivity and output matrix \texttt{Nodes2Edge} with dimension $\texttt{nC} \times \texttt{nC}$ (see Listing \ref{lst:VecNumEdges}).
\begin{equation}
\texttt{Nodes2Edge}(i,j)= \begin{cases} p ~ ~~ \text{if edge}~p~\text{between the nodes}~i~\text{and}~j; \\
0 ~~~\text{if there is no edge between the nodes}~i~\text{and}~ j.
\end{cases}
\end{equation}
The matrix \texttt{Nodes2Edge} is symmetric since it stores the indices of undirected edges.

\begin{lstlisting}[caption= $\mathbf{Vectorized ~NumEdges.m}$, label=lst:VecNumEdges, escapechar=|]
function [nEd,Nodes2Edge,U]=NumEdges(Coord,Tetra)
nC=size(Coord,1); |\label{line:ncNumEdges}|% Number of nodes/vertices in triangulation
I=Tetra(:,1); J=Tetra(:,2); K=Tetra(:,3); L=Tetra(:,4);% Collecting the nodes of tetrahedrons|\label{line:IJKLNumEdges}|
e1=[I J];e2=[I K];e3=[I L];e4=[J K];e5=[J L];e6=[K L]; % edges of all tetrahedrons |\label{line:e1e2e3e4e5e6NumEdges}|
r=reshape([e1(:) e2(:) e3(:) e4(:) e5(:) e6(:)]',[],2); |\label{line:rNumEdges}| % This matrix collects the initial and end node of edges
U=unique(sort(r,2),'rows','stable'); |\label{line:UNumEdges}|% removing repeated edges
nEd=length(U);  |\label{line:TotalEdgeNumEdges}| % total number of edges
S=sparse(U(:,1),U(:,2),1:nEd,nc,nc); |\label{line:SNumEdges}| % marking global edge number to edges
Nodes2Edge=S+S'; |\label{line:Nodes2Edge}|
end
\end{lstlisting}

\begin{example}
The full version of the sparse matrix \texttt{Nodes2Edge} for \figurename~\ref{fig:3DInitialMesh} is:
\[
\texttt{Nodes2Edge} =
\begin{bmatrix}
0 & \circnum{2} & \circnum{1} & \circnum{3} & 0 & 0 & 0 & 0 \\
2 & 0 & 4 & 6 & 11 & 10 & 0 & 18 \\
1 & 4 & 0 & 5 & 15 & 0 & 13 & 16 \\
3 & 6 & 5 & 0 & 8 & 7 & 14 & 0 \\
0 & 11 & 15 & 8 & 0 & 9 & 12 & 17 \\
0 & 10 & 0 & 7 & 9 & 0 & 0 & 0 \\
0 & 0 & 13 & 14 & 12 & 0 & 0 & 0 \\
0 & 18 & 16 & 0 & 17 & 0 & 0 & 0
\end{bmatrix}_{8 \times 8}
\]

\vspace{1em}

The first row of \texttt{Nodes2Edge} indicates:
\begin{itemize}
    \item The undirected edge $1 \leftrightarrow 2$ has number \circnum{2},
    \item The undirected edge $1 \leftrightarrow 3$ has number \circnum{1},
    \item The undirected edge $1 \leftrightarrow 4$ has number \circnum{3}.
\end{itemize}
\end{example}

\subsection{Edges-to-Tetrahedron connectivity}
The function \textbf{EdgeNum2Tetra} establishes the connectivity between edges and tetrahedra in a mesh (see Listing \ref{lst:VecEdge2NumTetra}). It returns a matrix \texttt{Edges2Tetra} with dimensions $\texttt{nEd} \times \texttt{nEd}$. The entries of \texttt{Edges2Tetra} indicate the tetrahedron numbers associated with each edge or zero if there is no tetrahedron connected to the edge.

We use the \figurename{~\ref{fig:3DTetrahedronsOrientation}} to illustrate the idea behind \textit{Edge-to-Tetrahedron} connectivity in the following examples,
\vspace{-0.3cm}
\begin{example}
If $e_1$ and $e_2$ are the edge numbers of the edges between the nodes $(a,b)$ and $(b,c)$, respectively. Then we have $\texttt{Edges2Tetra}(e_1,e_2)=T_{+}$ and $\texttt{Edges2Tetra}(e_2,e_1)=T_{-}$.
\end{example}
\begin{example}
If $e_1$ and $e_2$ are the edge numbers of the edges between the nodes $(a,c)$ and $(a,e)$, respectively. Then we have $\texttt{Edges2Tetra}(e_1,e_2)=T_{-}$ and $\texttt{Edges2Tetra}(e_2,e_1)=0$.
\end{example}

\begin{lstlisting}[caption= $\mathbf{Vectorized~ Edge2NumTetra.m}$, label=lst:VecEdge2NumTetra, escapechar=|]
function [Edges2Tetra]=EdgeNum2Tetra(Nodes2Edge,nEd,Tetra)
nE=size(Tetra,1);  % Number of tetrahedrons
I=Tetra(:,1); J=Tetra(:,2); K=Tetra(:,3); L=Tetra(:,4); % Collecting the nodes of tetrahedrons |\label{line:IJKLEdgeNum2Tetra}|
S=reshape([1:nE,1:nE,1:nE]',3*nE,1); |\label{line:SEdgeNum2Tetra}|
[ex,ey]=e2n(Nodes2Edge,I,J,K);|\label{line:StartEdges2TetraFirst}| % e1 is i-->j, e2 is j-->k and e3 is k-->i are edges first face 
Edges2Tetra=sparse(ex,ey,S,nEd,nEd); |\label{line:Edges2TetraFirst}|
[ex2,ey2]=e2n(Nodes2Edge,I,L,J);|\label{line:StartEdges2TetraSecond}| % e1 is i-->l, e2 is l-->j and e3 is j-->i are edges second face
Edges2Tetra=Edges2Tetra+sparse(ex2,ey2,S,nEd,nEd); |\label{line:Edges2TetraSecond}|
[ex3,ey3]=e2n(Nodes2Edge,I,K,L);|\label{line:StartEdges2TetraThird}| % e1 is i-->k, e2 is k-->l and e3 is l-->i are edges third face
Edges2Tetra=Edges2Tetra+sparse(ex3,ey3,S,nEd,nEd); |\label{line:Edges2TetraThird}|
[ex4,ey4]=e2n(Nodes2Edge,J,L,K);|\label{line:StartEdges2TetraFourth}| % e1 is j-->l, e2 is l-->k and e3 is k-->j are edges fourth face
Edges2Tetra=Edges2Tetra+sparse(ex4,ey4,S,nEd,nEd); |\label{line:Edges2TetraFourth}|
    function [ex,ey]=e2n(N2Ed,i,j,k)
        sz=size(N2Ed);
        ex=[N2Ed(sub2ind(sz,j,k)); N2Ed(sub2ind(sz,k,i)); N2Ed(sub2ind(sz,i,j))];
        ey=[N2Ed(sub2ind(sz,i,j)); N2Ed(sub2ind(sz,j,k)); N2Ed(sub2ind(sz,k,i))];    
    end
end
\end{lstlisting}

\begin{itemize}
    \item Lines \ref{line:StartEdges2TetraFirst}-\ref{line:Edges2TetraFirst}: Creating edges of the first oriented face of the tetrahedrons:  

    \begin{itemize}
    
    \item The vectors \texttt{ex} and \texttt{ey} represent the node indices of the three edges in the first oriented face of all tetrahedrons.
    
    \item 
    The \textit{Edge-to-Tetrahedron} connectivity over the first oriented face is stored in \texttt{Edges2Tetra} matrix.
    \end{itemize}
    \item Lines \ref{line:StartEdges2TetraSecond}-\ref{line:Edges2TetraFourth}: The code follows a similar pattern for the second, third, and fourth-oriented faces, but with different combinations of nodes and edge matrices. Finally, update the \texttt{Edges2Tetra} matrix.
    
\end{itemize}
\remark{For an $m^{\text{th}}$ tetrahedron, denoted $T_{m}=[a~b~c~d]$, consider one of its oriented face, say $F=[a~b~c]$. Let the edges $e_1,e_2$ and $e_3$ correspond to the pairs $(a,b),(b,c)$ and $(c,a)$ respectively. Then the $\texttt{Edges2Tetra}(e_2,e_1)=\texttt{Edges2Tetra}(e_3,e_2)=\texttt{Edges2Tetra}(e_1,e_3)=m$, but the  $\texttt{Edges2Tetra}(e_1,e_2)=\texttt{Edges2Tetra}(e_2,e_3)= \texttt{Edges2Tetra}(e_3,e_1) \neq m$ since these corresponds to the oriented face $F=[c~b~a]$, which is in opposite orientation with the tetrahedron $T_m$. Now, we can state that the 
\texttt{Edges2Tetra} is not a symmetric matrix.
}

\subsection{Faces-to-Tetrahedron connectivity}
The \textbf{FaceUP} function generates a matrix \texttt{Faces} with a dimension of $\texttt{nF} \times 3$ (see Listing \ref{lst:VecFace}). Let tetrahedron $T$ be denoted by the 4-tuple $[a~ b~ c~ d]$, then we define $m$-th row of the matrix \texttt{Faces} by,
\begin{align}
\texttt{Faces}(m,:)=[x~ y~ z]~~~~~~~~~~&\text{if}~\overrightarrow{yz}\times \overrightarrow{yx} ~\text{is an outward normal to  the}\\
& \text{tetrahedron}~T~\text{for}~x,y,z \in \{a,b,c,d\}. \nonumber
\end{align}
 If $[x~ y~ z]$ is an updated face, then all 3-tuples generated by $\{x,y,z\}$ without repetitions will not be updated. Thus, each row entry in \texttt{Faces} is a 3-tuple representing a unique orientation for a face in $\mathcal{F}_h$.

\begin{lstlisting}[caption= $\mathbf{Vectorized ~FaceUP.m}$, label=lst:VecFace, escapechar=|]
function [Faces]= FacesUp(Tetra, Edges2Tetra, Nodes2Edge)
nE=size(Tetra,1); |\label{line:ne}| % Total number of tetrahedrons in triangulation
I=Tetra(:,1); J=Tetra(:,2); K=Tetra(:,3); L=Tetra(:,4);%Collecting the nodes of tetrahedrons |\label{line:IJKLFace}|
I3=reshape(1:numel(Tetra),4,nE)'; % used to index the faces of the tetrahedrons|\label{line:I3Face}|
Q=sparse(4*nE,3); |\label{line:QFace}|
ind1=faceindex(Edges2Tetra,Nodes2Edge,I,J,K,[1:nE]'); |\label{line:ind1Face}|
Q(I3(ind1,1),:)=[Tetra(ind1,1) Tetra(ind1,2) Tetra(ind1,3)];%collecting 1st faces of all tetrahedrons  |\label{line:firstQupdateFace}| 
ind2=faceindex(Edges2Tetra,Nodes2Edge,I,K,L,[1:nE]');
Q(I3(ind2,2),:)=[Tetra(ind2,1) Tetra(ind2,3) Tetra(ind2,4)];%collecting 2nd faces of all tetrahedrons
ind3=faceindex(Edges2Tetra,Nodes2Edge,I,L,J,[1:nE]');
Q(I3(ind3,3),:)=[Tetra(ind3,1) Tetra(ind3,4) Tetra(ind3,2)];%collecting 3rd faces of all tetrahedrons
ind4=faceindex(Edges2Tetra,Nodes2Edge,L,K,J,[1:nE]');
Q(I3(ind4,4),:)=[Tetra(ind4,4) Tetra(ind4,3) Tetra(ind4,2)];%collecting 4th faces of all tetrahedrons |\label{line:EndQFace}|
Faces=Q(find(any(Q, 2)),:); |\label{line:Faces_vec}| 
    function ind=faceindex(E2T,N2Ed,i,j,k,p)
        sz=size(N2Ed); se=size(E2T);
        AdjTetra=E2T(sub2ind(se,N2Ed(sub2ind(sz,i,j)), N2Ed(sub2ind(sz,j,k))));
        ind=sort(nonzeros([find(AdjTetra==0);find(p<AdjTetra)]));
    end
end
\end{lstlisting}

\begin{itemize}
    \item Line \ref{line:QFace}-\ref{line:EndQFace}: We initialize a sparse matrix \texttt{Q} to store the faces data by using adjacent tetrahedra indices for each face. Use the \textbf{faceindex} function to identify these indices, compute the adjacent tetrahedra indices for each face in the tetrahedra, and update \texttt{Q} with the node indices of the adjacent tetrahedra for each face.
\item Line \ref{line:Faces_vec}: This line collects the rows of \texttt{Q} that contain at least one nonzero element using the \textbf{any} function. The resulting matrix \texttt{Faces} contains the node indices of the tetrahedrons forming the faces of the triangulation.
\end{itemize}

\begin{example}
    The \figurename{~\ref{fig:FaceUP_Faces_Matrix}} is the geometrical representation of the matrix \texttt{Faces} for 
    initial mesh data (see Listing \ref{lst:InitMeshData}) with dimension $16 \times 3$.
    The first four rows (red dotted) are the faces of the first tetrahedron $[1~2~3~5]$, rows $5^{\text{th}}~\text{to}~8^{\text{th}}$ (green dotted), $9^{\text{th}}~\text{to}~12^{\text{th}}$ (blue dotted), $13^{\text{th}}~\text{to}~16^{\text{th}}$ (cyan dotted) rows represent the faces of second $[5~8~6~2]$, third $[7~3~8~5]$ and fourth $[4~8~3~2]$ tetrahedron, respectively. 
    The fifth tetrahedron $[8~3~2~5]$ is the interior element of a mesh, so the faces $[8~3~2],~[8~2~5],~[8~5~3],~\text{and}~[5~2~3]$ formed by this tetrahedron are excluded in the \texttt{Faces} matrix because their tuples are already repeated in $16^{\text{th}},~7^{\text{th}},~12^{\text{th}},~\text{and}~4^{\text{th}}$ rows, respectively.
    \begin{figure}[H]             
\centering
\vspace{-0.25cm}
\hspace{1.0cm}
\includegraphics[height=7cm, width=7cm]{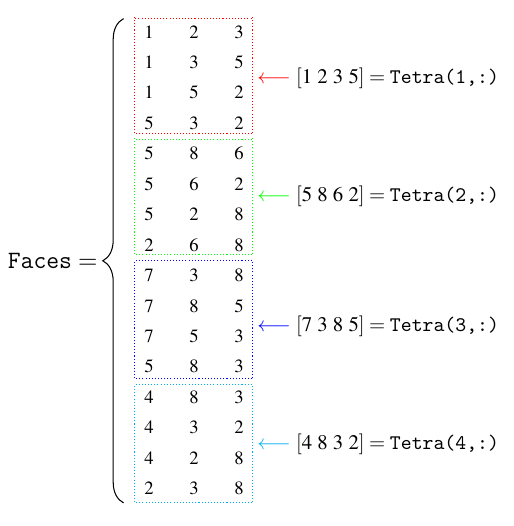}
\caption{\footnotesize Construction of matrix \texttt{Faces} for the \figurename{~\ref{fig:3DInitialMesh}}.}
 \label{fig:FaceUP_Faces_Matrix}
 \end{figure}
\end{example}

\subsection{RedRefinement}

\begin{lstlisting}[caption= $\mathbf{Vectorized~ RedRefine.m}$, label=lst:VecRedRefine, escapechar=|]
function [Tetra,Coord,A,Db,Nb]=RedRefine(Coord,Tetra,Db,Nb,U)
nC=size(Coord,1);  % Number of nodes
nE=size(Tetra,1); % Number of tetrahedrons
I=Tetra(:,1); J=Tetra(:,2); K=Tetra(:,3); L=Tetra(:,4);
M1=(Coord(I,:)+Coord(J,:))./2;  M2=(Coord(I,:)+Coord(K,:))./2;  M3=(Coord(I,:)+Coord(L,:))./2; |\label{line:M1M2M3RedRefine}| 
M4=(Coord(J,:)+Coord(K,:))./2; M5=(Coord(J,:)+Coord(L,:))./2; M6=(Coord(K,:)+Coord(L,:))./2;  |\label{line:M4M5M6RedRefine}| 
Ct=(Coord(I,:)+Coord(J,:)+Coord(K,:)+Coord(L,:))./4;  |\label{line:centridRedRefine}| % centroid 
r1=reshape([Ct(:) M1(:) M2(:) M3(:) M4(:) M5(:) M6(:)]', [], 3); |\label{line:r1RedRefine}|  
Coord=[Coord; unique(r1,'stable','rows')]; |\label{line:CoordRedRefine}|  % update Coord with coordinates of new nodes
C = find(ismember(Coord, Ct, 'rows')); |\label{line:CRedRefine}| 
R=setdiff((nC+1):length(Coord),C); % indices of the newly formed nodes|\label{line:RRedRefine}| 
S=sparse(U(:,1),U(:,2),R,nC,nC); |\label{line:SRedRefine}| 
A=S+S'; w=size(A); |\label{line:ARedRefine}|  % This matrix collects the number of newly formed nodes
%% Forms 4 refined tetrahedrons at 4 newly formed nodes with the midpoint of edges
Tetra(1:nE,:)=[I A(sub2ind(w,I,J)) A(sub2ind(w,I,K)) A(sub2ind(w,I,L))] ; |\label{line:TetraStartRedRefineMidpoint}| 
Tetra(nE+1:11:12*nE,:)=[J A(sub2ind(w,K,J)) A(sub2ind(w,I,J)) A(sub2ind(w,J,L))];  
Tetra(nE+2:11:12*nE,:)=[K A(sub2ind(w,I,K)) A(sub2ind(w,J,K)) A(sub2ind(w,K,L))]; 
Tetra(nE+3:11:12*nE,:)=[L A(sub2ind(w,K,L)) A(sub2ind(w,L,J)) A(sub2ind(w,I,L))]; |\label{line:TetralastRedRefineMidpoint}|
%% These are the rest 8 tetrahedrons contributed by the octahedron |\label{line:TetraStartRedRefineCentroid}|
Tetra(nE+4:11:12*nE,:)=[A(sub2ind(w,I,J)) A(sub2ind(w,I,K)) A(sub2ind(w,I,L)) C]; 
Tetra(nE+5:11:12*nE,:)=[A(sub2ind(w,K,J)) A(sub2ind(w,I,J)) A(sub2ind(w,J,L)) C]; 
Tetra(nE+6:11:12*nE,:)=[A(sub2ind(w,I,K)) A(sub2ind(w,J,K)) A(sub2ind(w,K,L)) C]; 
Tetra(nE+7:11:12*nE,:)=[A(sub2ind(w,K,L)) A(sub2ind(w,L,J)) A(sub2ind(w,I,L)) C]; 
Tetra(nE+8:11:12*nE,:)=[A(sub2ind(w,I,J)) A(sub2ind(w,K,J)) A(sub2ind(w,I,K)) C]; 
Tetra(nE+9:11:12*nE,:)=[A(sub2ind(w,I,K)) A(sub2ind(w,L,K)) A(sub2ind(w,I,L)) C]; 
Tetra(nE+10:11:12*nE,:)=[A(sub2ind(w,I,J)) A(sub2ind(w,I,L)) A(sub2ind(w,J,L)) C]; 
Tetra(nE+11:11:12*nE,:)=[A(sub2ind(w,K,L)) A(sub2ind(w,K,J)) A(sub2ind(w,L,J)) C];  |\label{line:TetralastRedRefineCentroid}| 
%% Boundary conditions
if(~isempty(Db)) Db=refineFaces(Db,A); end |\label{line:StartDiricheltBoundary}|  |\label{line:LastDiricheltBoundary}|
if(~isempty(Nb)) Nb=refineFaces(Nb,A); end  |\label{line:StartNeumanntBoundary}| |\label{line:LastNeumannBoundary}|
    function F2n=refineFaces(F2n,A)
        n1=F2n(:,1); n2=F2n(:,2); n3=F2n(:,3); 
        nF=size(F2n,1); w=size(A);
        F2n(1:nF,:)=[n1 A(sub2ind(w,n1,n2)) A(sub2ind(w,n1,n3))];
        F2n(nF+1:3:4*nF,:)=[n2 A(sub2ind(w,n2,n3)) A(sub2ind(w,n1,n2))]; 
        F2n(nF+2:3:4*nF,:)=[n3 A(sub2ind(w,n1,n3)) A(sub2ind(w,n2,n3))]; 
        F2n(nF+3:3:4*nF,:)=[A(sub2ind(w,n1,n2)) A(sub2ind(w,n2,n3)) A(sub2ind(w,n1,n3))];
    end |\label{line:endRefineFaces}|
end
\end{lstlisting}
Here, we discuss the refinement of tetrahedrons in Listing \ref{lst:VecRedRefine}, and the \figurename{~\ref{fig:CoarshRefinedOctahedron}} can be reproduced by the script
$$\texttt{fig\_CoarshRefinedOctahedron.m}$$
which illustrates how a single tetrahedron (see \figurename{~\ref{fig:sub1}}) is refined into 12 tetrahedrons. First, nodes are created at the edge's midpoints and the tetrahedron's centroid, resulting in a total of seven new nodes formed in a single tetrahedron (see lines \ref{line:M1M2M3RedRefine}-\ref{line:CoordRedRefine}). Each node in the coarser tetrahedron, by using the nodes at the mid-points of the edges that share the node, forms a refined tetrahedron, leading to 4 refined tetrahedrons (see lines \ref{line:TetraStartRedRefineMidpoint}-\ref{line:TetralastRedRefineMidpoint}). If we chop off these newly formed tetrahedrons, we will be left with an octahedron, which can be seen on the \figurename{~\ref{fig:sub3}. Now, join all 8 faces present on the boundary of the octahedron with the node at the centroid of the tetrahedron, which will lead to 8 more refined tetrahedrons (see lines \ref{line:TetraStartRedRefineCentroid}-\ref{line:TetralastRedRefineCentroid}). This way, we will refine each coarser tetrahedron into 12 refined tetrahedrons (see \figurename{~\ref{fig:sub2}}). This construction ensures there won't be any possibility of hanging nodes in the refinement since new nodes are created at the edge's midpoint or the centroid of the tetrahedron. The same idea can be seen in the function \textbf{RedRefine.m}.
For the boundary faces ($\mathcal{F}^h_{D} \cup \mathcal{F}^h_{N}$), new nodes are created at the centroid of the boundary faces (see lines \ref{line:StartDiricheltBoundary}-\ref{line:endRefineFaces}). Update boundary data with the new nodes.


\begin{figure}[H]
    \centering
    \subfigure[Coarser tetrahedron]{
        \includegraphics[width=0.305\textwidth]{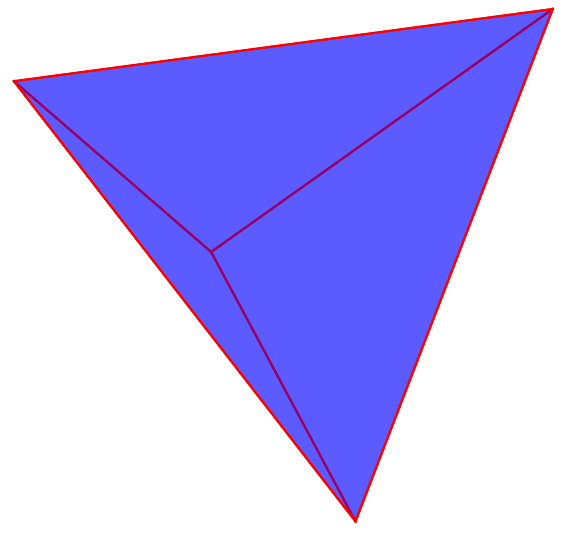}
        \label{fig:sub1}
    }
    \hfill
    \subfigure[Refined tetrahedron]{
        \includegraphics[width=0.305\textwidth]{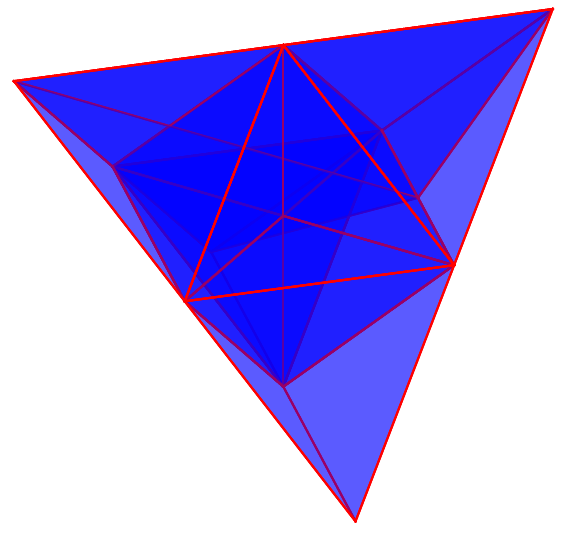}
        \label{fig:sub2}
    }
    \subfigure[Refined Octahedron]{
        \includegraphics[width=0.305\textwidth]{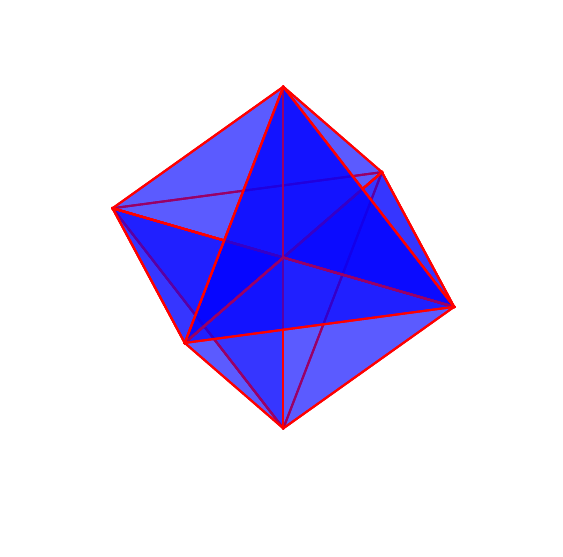}
        \label{fig:sub3}
    }
    \caption{Red refinement of a tetrahedron: subdivision into twelve smaller tetrahedrons.}
    \label{fig:CoarshRefinedOctahedron}
\end{figure}

\begin{example}
    Consider the coarser tetrahedron $\texttt{T}=[1~3~2~4]$ (cf. \figurename{~\ref{fig:sub1}}) with vertices at coordinates $\texttt{C}=[0~0~0; 1~0~0; 0~1~0; 0~0~1]$. We then refine $\texttt{T}$ into 12 smaller tetrahedrons (cf. \figurename{~\ref{fig:sub2}}) by introducing 7 new nodes in $\texttt{C}$ and updating data's accordingly,
\begin{multicols}{2}
\hspace{3cm}
$\texttt{T}=\begin{bmatrix}
1 & 6 & 7 & 8 \\
3 & 9 & 6 & 10 \\
2 & 7 & 9 & 11 \\
4 & 11 & 10 & 8 \\
6 & 7 & 8 & 5 \\
9 & 6 & 10 & 5\\
7 & 9 & 11 & 5 \\
11 & 10 & 8 & 5 \\
6 & 9 & 7 & 5 \\
7 & 11 & 8 & 5 \\
6 & 8 & 10 & 5 \\
11 & 9 & 10 & 5 \\
\end{bmatrix}_{12 \times 4}$,
\hspace{3cm}
$\texttt{C}=
\begin{bmatrix}
0&0&0 \\ 
1&0&0 \\
0&1&0\\
0&0&1\\\vspace{0.15cm}
\frac{1}{4}&\frac{1}{4}&\frac{1}{4} \\ \vspace{0.15cm}
0 & \frac{1}{2} & 0 \\ \vspace{0.15cm}
\frac{1}{2}  & 0 & 0 \\ \vspace{0.15cm}
0& 0 &  \frac{1}{2}  \\ \vspace{0.15cm}
\frac{1}{2} &\frac{1}{2}  & 0 \\ \vspace{0.15cm}
0 & \frac{1}{2}  & \frac{1}{2}  \\ \vspace{0.15cm}
\frac{1}{2}  & 0 & \frac{1}{2}  \\
\end{bmatrix}_{11 \times 3}$.
\end{multicols}

\end{example}

\subsection{Run-time of \textbf{FaceUp.m}, \textbf{NumEdges.m}, \textbf{EdgeNum2Tetra.m}, and \textbf{RedRefine.m}}

Table~\ref{tab:RefinementPythonMATLAB} presents the run-time performance of the vectorized implementations of the functions used for uniform refinement of the initial triangulation. For illustration purposes, the meshes obtained after two and three levels of refinement are depicted in Figure~\ref{3Drefine}. After five levels of refinement, the resulting mesh consists of 1,244,160 tetrahedra.

\begin{table}[h!]
  \begin{center}
  \centering
  \caption{Vectorized code run-time of the \textbf{FaceUP}, \textbf{NumEdges}, \textbf{EdgeNum2Tetra}, and \textbf{RedRefine} function.}
  \footnotesize\addtolength{\tabcolsep}{0pt}
    \label{tab:RefinementPythonMATLAB}
    \begin{tabular}{|c|r|r|r|r|c|c|c|c|}  
      \hline 
    Level & \texttt{nE} & \texttt{nC} &\texttt{nEd} & \texttt{nF} & \multicolumn{4}{c|}{Computational time (in sec.)  for}   \\ \cline{6-9}
     & & & & & \textbf{FaceUp.m}  &  \textbf{NumEdges.m}  & \textbf{EdgeNum2Tetra.m}  & \textbf{RedRefine.m}  \\
      \hline 
    1 & 60 & 31 & 114       & 144& 0.0005 & 0.0013  & 0.0010  & 0.0017 \\
    2 & 720   & 205 & 1,020     & 1,536 &   0.0025 & 0.0030 & 0.0023   & 0.0108\\
     3 & 8,640 & 1945 & 10,968    & 17,664 & 0.0716  & 0.0091 & 0.0162 & 0.0168\\
    4 & 103,680 & 21,553 & 126,768   & 208,896 &   0.2320  &  0.0912  &  0.1601 & 0.0480\\
    5 & 1,244,160 & 252,001 & 1,502,304 & 2,494,464&  3.7920  &  1.0248 &  2.7785 &  0.5407\\        
      \hline 
    \end{tabular} 
  \end{center}
\end{table}

\begin{figure}
 \centering 
   \includegraphics[scale=0.65]{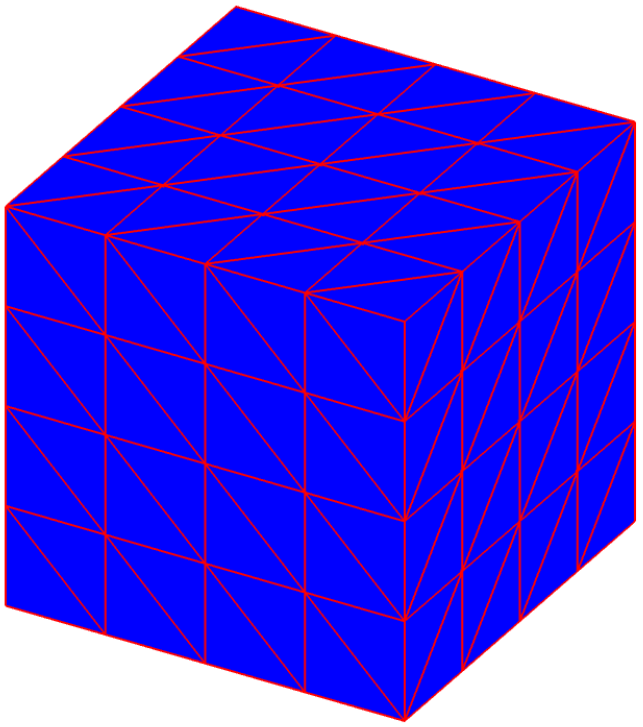}
   \includegraphics[scale=0.65]{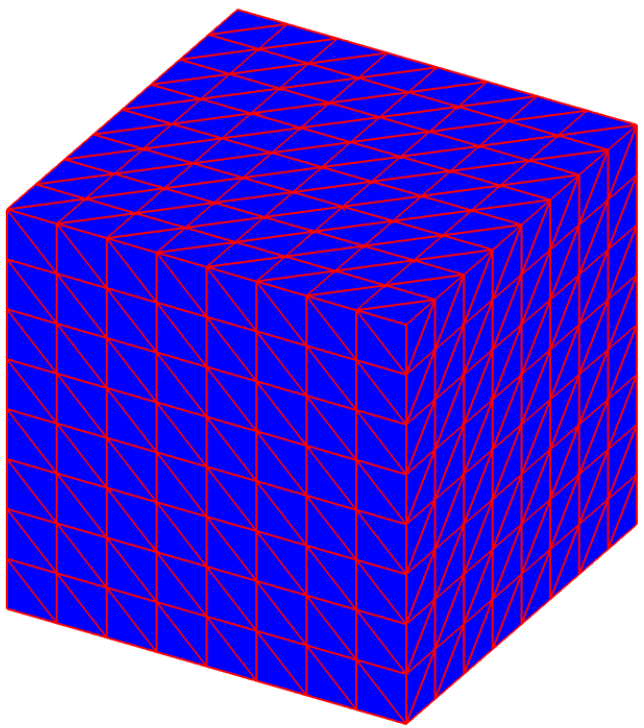}
  \caption{Level 2 (left) and level 3 (right) refined mesh.}
  \label{3Drefine}
\end{figure}

\section{Primal hybrid implementation of second order elliptic problem}\label{sub:primalhybrid}
We consider the model problem: find $u:\Omega\to\mathbb{R}$ such that
\begin{equation}\label{eq1}
\begin{cases}
\begin{alignedat}{2}
-\nabla\!\cdot(\nabla u) + u & = f  &\quad & \text{in } \Omega,\\[4pt]
u                           & = u_D &      & \text{on } \Gamma_D,\\[4pt]
\nabla u \cdot \vec{\nu}    & = g   &      & \text{on } \Gamma_N,
\end{alignedat}
\end{cases}
\end{equation}
where $f$, $u_D$, and $g$ are sufficiently smooth and $\vec{\nu}$ is the unit outward normal on $\partial\Omega$.
We use the standard Sobolev spaces $H^{1}(\Omega)$ and 
$H(\texttt{div},\Omega)$ defined in the domain $\Omega$.
Let $H^1_{D}(\Omega)=\{v \in H^1(\Omega): v_{|_{\partial D} }=u_D\}$.
 We denote by $H^{1/2}(\Gamma)$ the space of traces ${v|}_{\Gamma}$
over $\Gamma$ of the functions $v\in H^{1}(\Omega)$. Let $H^{-1/2}(\Gamma)$ denote the dual space of $H^{1/2}(\Gamma)$ and $\langle\cdot,\cdot\rangle_{\Gamma}$ denote their pairing.
We consider the broken Sobolev space
\begin{equation*}
\mathbb{Y}=\{v\in L^{2}(\Omega): {v|}_{T}\in H^{1}(T), \quad \forall T\in\mathcal{T}_{h}\}=\displaystyle\prod_{T\in\mathcal{T}_{h}}H^{1}(T)
\end{equation*}
equipped with the mesh-dependent norm \begin{equation}
\|v\|_\mathbb{Y}=\bigg(\displaystyle\sum_{T\in\mathcal{T}_h}\|\grad( v)\|^2_{L^2(T)}+h^{-2}\|v\|^{2}_{L^2(T)}\bigg)^{\frac{1}{2}}.\end{equation}
We define the Lagrange multiplier space as
\begin{align*}
\mathbb{S}=\bigg\{\chi\in\prod_{T\in\mathcal{T}_{h}}H^{-1/2}(\partial T/\Gamma_N):~\exists~\Vec{q}\in H(\texttt{div},\Omega)
~\text{such that},
~\Vec{q}\cdot\Vec{\nu}_{T}=\chi ~\text{on} ~\partial T/\Gamma_N, \forall T\in\mathcal{T}_{h}\bigg\}
\end{align*}
equipped with the norm  
(cf. \cite{raviart}) 
\begin{equation}\label{eMnorm}
   \|\chi \|_h =\bigg( \sum_{T \in \mathcal{T}_h} h \|\chi \|^2_{L^2(\partial T/\Gamma_N)}  \bigg)^{\frac{1}{2}} .
\end{equation}
 
Primal hybrid formulation of  \eqref{eq1} is to seek a pair of solutions  
$(u,\kappa)\in \mathbb{Y}\times \mathbb{S}$ such that
\begin{alignat}{3}
&\sum_{T\in\mathcal{T}_{h}} \int_{T} \bigl( \nabla u \cdot \nabla v + u v \bigr)\dx
  -\sum_{F\in \mathcal{F}_{h}\setminus \mathcal{F}^h_{N}} \langle \kappa, v \rangle_{F} 
&\;=\;& \mathbb{F}(v), 
&\qquad & \forall v \in \mathbb{Y}, \label{eq1hyb1} \\[1mm]
-&\sum_{F\in \mathcal{F}_{h}\setminus \mathcal{F}^h_{N}} \langle \chi, u \rangle_{F} 
&\;=\;& -\sum_{F\in \mathcal{F}_{h}\setminus \mathcal{F}^h_{N}} \langle \chi, u_D \rangle_{F}, 
& & \forall \chi \in \mathbb{S}. \nonumber
\end{alignat}
Here, the Lagrange multiplier associated with the constraint $u\in H^{1}_{D}(\Omega)$ reads
\begin{equation}
\kappa=\grad (u)\cdot\Vec{\nu}_T \quad \mbox{on}~ \partial T/\Gamma_N, \quad \forall T\in\mathcal{T}_{h}    
\end{equation}
and 
$$\mathbb{F}(v)=\displaystyle \sum_{T\in\mathcal{T}_h}\int_Tfv\dx+\int_{\Gamma_N} gv\dgamma.$$
\subsection{Discrete problem and algebraic formulation}

Let $P_k$ be the space of polynomials of degree $\leq k$. We define the finite-dimensional spaces $\mathbb{Y}_h\subset \mathbb{Y}$ and $\mathbb{S}_h\subset \mathbb{S}$ as 
\begin{equation} \label{eq:Primalspace}
\displaystyle  \mathbb{Y}_h=\prod_{T\in\mathcal{T}_h}P_1(T)   
\end{equation}
and
\begin{align}\label{eq:HybridSpace}
\mathbb{S}_{h}=\bigg\{\chi_{h}\in\prod_{\partial T\in \Gamma/\Gamma_{N}}P_0(\partial T):{\chi_{h}|}_{\partial T_{+}}
+{\chi_{h}|}_{\partial T_{-}}=0~\text{on}~T_{+}\cap T_{-},~
\text{for every adjacent pair}~T_{\pm}\in\mathcal{T}_{h}\bigg\}.
\end{align}

We seek a pair $(u_h,\kappa_h)\in \mathbb{Y}_h\times \mathbb{S}_h$ such that

\begin{alignat}{3}
&\sum_{T\in\mathcal{T}_h} \int_{T} \bigl( \nabla u_h \cdot \nabla v_h + u_h v_h \bigr) \dx
  - \sum_{F\in \mathcal{F}_{h} \setminus \mathcal{F}^h_{N}} \int_{F} \kappa_h [\![ v_h ]\!]_F \, \dgamma 
&\;=\;& \mathbb{F}(v_h), 
&\quad & \forall v_h \in \mathbb{Y}_h, \label{ph21} \\[1mm]
&-\sum_{F\in \mathcal{F}_{h} \setminus \mathcal{F}^h_{N}} \int_{F} \chi_h [\![ u_h ]\!]_F \, \dgamma
&\;=\;& -\sum_{F\in \mathcal{F}_{h} \setminus \mathcal{F}^h_{N}} \int_{F} \chi_h u_{Dh} \dgamma, 
& & \forall \chi_h \in \mathbb{S}_h. \label{ph22}
\end{alignat}
Note that problem \eqref{ph21}-\eqref{ph22} has a unique pair of solutions $(u_h,\kappa_h)\in \mathbb{Y}_h\times \mathbb{S}_h$ (cf. \cite{raviart}).

 With $N=dim(\mathbb{Y}_h)=4*\texttt{nE}$, let $\{\phi_1,\cdots,\phi_N\}$ be the basis for $\mathbb{Y}_h$ and $U=[\alpha_1,\cdots,\alpha_N]^\prime$ be the components of $$u_h=\displaystyle\sum^N_{i=1}\alpha_i\phi_i.$$

Let $F_1, F_2, F_3, F_4$ be the faces of the tetrahedron $T$,
and let $\Vec{\nu}_{F_k}$ denote the unit normal vector of $F_k$ chosen with a global fixed orientation while $\Vec{\nu}_k$ denotes the outer unit normal of $T$ along $F_k$.
We define the basis functions $\psi_F$ of $\mathbb{S}_h$, where $F\in\mathcal{F}_h/\mathcal{F}^h_{N}$ as follows:
\begin{equation}\label{eq:sigma}
\psi_{F_k}(\textbf{x}) = \sigma_k ~~~ \mathrm{for}~ k=1,2,3,4~ \mathrm{and}~\textbf{x}\in T,
\end{equation}
where $\sigma_k = \Vec{\nu}_{F_k}.\Vec{\nu}_{k}$ is $+1$ if $\Vec{\nu}_{F_k}$ points outward and otherwise $-1$.

Globally, let for any face $F=T_+\cap T_-\in\mathcal{F}^h_\Omega$
\begin{align*}
\psi_F(\textbf{x})=\begin{cases}
& 1~~~~~\text{for}~~\textbf{x}\in T_+,\\
& -1~~~\text{for}~~\textbf{x}\in T_-,\\
&0~~~~~~\text{elsewhere};
\end{cases}
\end{align*}
and for $F\in \mathcal{F}^h_D$
\begin{align*}
\psi_F(\textbf{x})=\begin{cases}
& 1~~\text{for}~~\textbf{x}\in F,\\
&0~~\text{elsewhere}.
\end{cases}
\end{align*}

 With $L=dim(\mathbb{S}_h)$ , let $\mathbb{S}_h=span\{\psi_l\}_1^L$ with $\psi_l=\psi_{F_l}$, where $l=1,\cdots, L$ is an enumeration of faces $\mathcal{F}^h_{\Omega}\cup\mathcal{F}^h_{D}=\{F_{1},\cdots,F_{L}\}$.
 Let $\Lambda=[\alpha_{N+1},\cdots,\alpha_{N+L}]^\prime$ be the components of $\kappa_h$ such that $$\kappa_h=\displaystyle\sum^{L}_{l=1}\alpha_{N+l}\psi_l.$$

 The problem \eqref{ph21}-\eqref{ph22} can be written in a linear system of equations for the unknowns $U$ and $\Lambda$ as
\begin{align}
&\sum^N_{i=1}\alpha_i\sum_{T\in\mathcal{T}_h}\int_T \left(\grad(\phi_j) \cdot\grad(\phi_i) + \phi_j \phi_i \right) \dx-\sum_{l=1}^{L}\alpha_{N+l}\int_{F_{l}}\psi_l[\![\phi_j]\!]_{F_l} \dgamma=\mathbb{F}(\phi_j),\label{ph23}\\
-&\sum_{i=1}^{N}\alpha_{i}\int_{F_l}\psi_l[\![\phi_i]\!]_{F_l} \dgamma\hspace{6.1cm}= - \int_{F_l}\psi_lu_{Dh} \dgamma,\label{ph24}
\end{align}

for $j=1,\cdots,N$ and $l=1,\cdots,L$.

\subsection{Local matrices}
For a tetrahedron element $T$, let $(x_1,y_1,z_1),~ (x_2,y_2,z_2),~(x_3,y_3,z_3),$ and $(x_4,y_4,z_4)$ be the vertices and $\phi_1,\phi_2,\phi_3,$ and $\phi_4$ be the corresponding basis functions.  

Then, the local stiffness and mass matrix,
\begin{equation}\label{eq:K(T)}
K(T) = \frac{|T|}{6} G G^{\mathrm{T}} ~~~~~~\mathrm{with}~~~~~G:= \begin{pmatrix}
         1 & 1 & 1 & 1\\
         x_1 & x_2 & x_3 & x_4\\
         y_1 & y_2 & y_3 & y_4\\
         z_1 & z_2 & z_3 & z_4\\
    \end{pmatrix}^{-1}   \begin{pmatrix}
         0 & 0 & 0\\
         1 & 0 & 0\\
         0 & 1 & 0\\
         0 & 0 & 1\\
    \end{pmatrix},
\end{equation}

\begin{equation}\label{eq:M(T)}
M(T) = \frac{|T|}{20} \begin{pmatrix}
         2 & 1 & 1 & 1\\
         1 & 2 & 1 & 1\\
         1 & 1 & 2 & 1\\
         1 & 1 & 1 & 2\\
    \end{pmatrix}.
\end{equation}

The local Lagrange multiplier matrix
\begin{align} \label{eq:C(T)}
C(T) &= \begin{pmatrix}
        \displaystyle \int_{F_{1}}\psi_1[\![\lambda_1]\!]_{F_1} \dgamma & \displaystyle \int_{F_{1}}\psi_1[\![\lambda_2]\!]_{F_1} \dgamma & \displaystyle \int_{F_{1}}\psi_1[\![\lambda_3]\!]_{F_1} \dgamma & \displaystyle \int_{F_{1}}\psi_1[\![\lambda_4]\!]_{F_1} \dgamma\\ 
        \displaystyle \int_{F_{2}}\psi_2[\![\lambda_1]\!]_{F_2} \dgamma & \displaystyle \int_{F_{2}}\psi_2[\![\lambda_2]\!]_{F_2} \dgamma & \displaystyle \int_{F_{2}}\psi_2[\![\lambda_3]\!]_{F_2} \dgamma & \displaystyle \int_{F_{2}}\psi_2[\![\lambda_4]\!]_{F_2} \dgamma\\
\displaystyle \int_{F_{3}}\psi_3[\![\lambda_1]\!]_{F_3} \dgamma & \displaystyle \int_{F_{3}}\psi_3[\![\lambda_2]\!]_{F_3} \dgamma & \displaystyle \int_{F_{3}}\psi_3[\![\lambda_3]\!]_{F_3} \dgamma & \displaystyle \int_{F_{3}}\psi_3[\![\lambda_4]\!]_{F_3} \dgamma\\
\displaystyle \int_{F_{4}}\psi_4[\![\lambda_1]\!]_{F_4} \dgamma & \displaystyle \int_{F_{4}}\psi_4[\![\lambda_2]\!]_{F_4} \dgamma & \displaystyle \int_{F_{4}}\psi_4[\![\lambda_3]\!]_{F_4} \dgamma &\displaystyle  \int_{F_{4}}\psi_4[\![\lambda_4]\!]_{F_4} \dgamma\\
\end{pmatrix} \nonumber \\
&=  \begin{pmatrix}
\frac{\sigma_1|F_1|}{3} & \frac{\sigma_1|F_1|}{3} & \frac{\sigma_1|F_1|}{3} & 0 \\
\frac{\sigma_2|F_2|}{3}  & 0 & \frac{\sigma_2|F_2|}{3}  & \frac{\sigma_2|F_2|}{3} \\
\frac{\sigma_3|F_3|}{3}  & \frac{\sigma_3|F_3|}{3} & 0  & \frac{\sigma_3|F_3|}{3}  \\
0 & \frac{\sigma_4|F_4|}{3} & \frac{\sigma_3|F_3|}{3} & \frac{\sigma_3|F_3|}{3}
\end{pmatrix} = \begin{pmatrix}
\sigma_1|F_1| r_{1} \\
\sigma_2|F_2| r_{2}  \\
\sigma_3|F_3| r_{3}   \\
\sigma_4|F_4| r_{4} 
\end{pmatrix}
\end{align}
where $r_{1}=\Big[\frac{1}{3}~\frac{1}{3}~\frac{1}{3}~0 \Big],~r_{2}=\Big[ \frac{1}{3}~0~\frac{1}{3}~\frac{1}{3} \Big],~r_{3}=\Big[ \frac{1}{3}~\frac{1}{3}~0~\frac{1}{3} \Big],~r_{4}=\Big[0~\frac{1}{3}~\frac{1}{3}~\frac{1}{3} \Big].$

\subsection{Assembly of global matrices}
The right-hand side part:
\begin{align*}
b_j&= \sum_{T \in \mathcal{T}_h} \int_{T} f \phi_{j} \dx 
\approx \displaystyle \sum_{T \in \mathcal{T}_h} \displaystyle \frac{|T|}{4} \displaystyle \sum_{i=1}^{4} f(Ct_{i}) \phi_{j}(Ct_{i}), \\
LN_j&=  \int_{\Gamma_{N}} g \phi_{j} \dgamma 
\approx   \displaystyle |F_i|  g(Ct_{i}) \phi_{j}(Ct_{i}), \\
 b_{D_j}&=\int_{F_l} \psi_l u_{Dh}  \dgamma 
\approx \displaystyle |F_l| \displaystyle u_{Dh}(Ct_l). 
\end{align*}

In vector-matrix form, \eqref{ph23}-\eqref{ph24} can be written as,
\begin{equation}\label{eq:VectorForm}
\left(
    \begin{array}{cc}
      K+M & -C^{\prime}  \\
      -C & 0  \\
    \end{array}
  \right)_{(N+L)\times (N+L)}  
  \quad  \left(\begin{array}{c}
      U   \\
      \Lambda  \\
    \end{array}
  \right)_{(N+L)\times 1} 
  = \left(
    \begin{array}{c}
      b+LN   \\
      -b_D  \\
    \end{array}
  \right)_{(N+L)\times 1}     
  \end{equation}
where

\begin{align*}
&K=\big(K_{ij}\big)_{N\times N}~, \quad K_{ij}=\displaystyle\sum_{T\in\mathcal{T}_h}\int_{T}\nabla \phi_i\cdot\nabla\phi_j~\dx, \\
&M=\big(M_{ij}\big)_{N\times N}~, \quad M_{ij}=\displaystyle\sum_{T\in\mathcal{T}_h}\int_{T} \phi_i\cdot\phi_j~\dx, \\
&C=(C_{N+l,i})_{L\times N}, \quad 
C_{N+l,i}=\displaystyle \int_{E_{l}}\psi_l[\![\phi_i]\!]_{E_l}\dgamma,\\
&b=(b_j)_{ N\times 1}, \quad LN=(LN_j)_{ N\times 1}, \quad b_D=(b_{D_l})_{L\times 1}.
\end{align*}

\subsection{StiffMassPH.m}
The function \textbf{StiffMassPH.m} (Listing \ref{lst:stiffPH}) implements the assembly of the stiffness matrix $K$, mass matrix $M$, and load vector $b$.
\begin{lstlisting}[caption= $\mathbf{StiffMassPH.m}$, label=lst:stiffPH, escapechar=|]
function [K,b,M,KT,MT,GI] = StiffMassPH(Tetra,Coord)
nE=size(Tetra,1); 
[vol,grad]=VolGrad(Coord,Tetra); |\label{line:volGrad}|
G1=reshape(grad(1,:,:),[nE,3]); G2=reshape(grad(2,:,:),[nE,3]); |\label{line:G1}|
G3=reshape(grad(3,:,:),[nE,3]); G4=reshape(grad(4,:,:),[nE,3]); |\label{line:G4}|
KT=[[dot(G1,G1,2),dot(G1,G2,2),dot(G1,G3,2),dot(G1,G4,2)].*vol; |\label{line:KTstart}|
    [dot(G2,G1,2),dot(G2,G2,2),dot(G2,G3,2),dot(G2,G4,2)].*vol;
    [dot(G3,G1,2),dot(G3,G2,2),dot(G3,G3,2),dot(G3,G4,2)].*vol;
    [dot(G4,G1,2),dot(G4,G2,2),dot(G4,G3,2),dot(G4,G4,2)].*vol]; |\label{line:KT}|
GI=reshape([1:4*nE; 1:4*nE; 1:4*nE; 1:4*nE],16,nE);  |\label{line:GI}|
I3=reshape(GI([1 5 9 13 2 6 10 14 3 7 11 15 4 8 12 16],:)',16*nE,1);   |\label{line:I3}|
J3=reshape(GI',16*nE,1); |\label{line:J3}|
K=sparse(I3,J3,KT,4*nE,4*nE); |\label{line:phK}|
%%%%%% Mass matrix %%%%%%%
MT=[[1/10,1/20,1/20,1/20].*vol;
   [1/20,1/10,1/20,1/20].*vol;
   [1/20,1/20,1/10,1/20].*vol;
   [1/20,1/20,1/20,1/10].*vol];
M=sparse(I3,J3,MT,4*nE,4*nE);  
%%%%%% Load vector %%%%%%%%%
P1=Coord(Tetra(:,1),:); P2=Coord(Tetra(:,2),:);
P3=Coord(Tetra(:,3),:); P4=Coord(Tetra(:,4),:);
Ct1=(P1+P2+P3)./3; Ct2=(P1+P3+P4)./3; Ct3=(P1+P4+P2)./3; Ct4=(P2+P4+P3)./3; % centroids of faces
f1=vol.*(f(Ct1)+f(Ct2)+f(Ct3))./12; f2=vol.*(f(Ct1)+f(Ct3)+f(Ct4))./12; |\label{line:f1}|
f3=vol.*(f(Ct1)+f(Ct4)+f(Ct2))./12; f4=vol.*(f(Ct2)+f(Ct4)+f(Ct3))./12; |\label{line:f4}| 
b=accumarray(reshape(1:4*nE,4*nE,1),reshape([f1 f2 f3 f4]',4*nE ,1),[4*nE 1]); |\label{line:phb}|
end
\end{lstlisting}

\begin{itemize}
   
    \item Line \ref{line:volGrad}: The function \textbf{VolGrad.m} computes the volume and gradients of the basis functions (see, \eqref{eq:K(T)}). 
    \item Lines \ref{line:KTstart}-\ref{line:KT}: The array
    \texttt{KT} stores the local stiffness matrix $K(T)$ (see \eqref{eq:K(T)})for each $T \in \mathcal{T}_h$.
    \item Lines \ref{line:GI}-\ref{line:I3}: We call the \texttt{GI} is global indices matrix with size $16\times \texttt{nE}$. The \texttt{I3} and \texttt{J3} are constructed based on the \texttt{GI} matrix. \newline
    Define, a row vector of size $\texttt{nE}\times 1$
    \begin{align*}
    H_{\alpha}&=\begin{bmatrix}
        4(1-1)+\alpha & 4(2-1)+\alpha & 4(3-1)+\alpha & \cdots & 4(\texttt{nE}-1)+\alpha
    \end{bmatrix},~\text{where}~\alpha=1,2,3,4 \\
    &=\begin{bmatrix}
        \alpha & 4+\alpha & 8+\alpha & \cdots & 4(\texttt{nE}-1)+\alpha
    \end{bmatrix},~\text{where}~\alpha=1,2,3,4.
    \end{align*}
    Then, using the definition of $H_{\alpha}$, we represent the entries of the matrix \texttt{GI} and the column vectors \texttt{I3} and \texttt{J3} as:

\[
\begin{tabular}{c@{\hspace{20pt}}@{}@{\hspace{20pt}} c @{\hspace{20pt}}@{}@{\hspace{20pt}}c@{\hspace{20pt}}} $\texttt{GI} = \begin{bmatrix}
    H_1 \\ H_1 \\ H_1 \\ H_1 \\ H_2 \\ H_2 \\ H_2 \\ H_2 \\ H_3 \\ H_3 \\ H_3 \\ H_3 \\ H_4 \\ H_4 \\ H_4 \\ H_4 
\end{bmatrix}_{16 \times \texttt{nE}}$ \hspace{-1cm}&
$\texttt{I3} = \begin{bmatrix}
    H^t_1 \\ H^t_2 \\ H^t_3 \\ H^t_4 \\ H^t_1 \\ H^t_2 \\ H^t_3 \\ H^t_4 \\ H^t_1 \\ H^t_2 \\ H^t_3 \\ H^t_4 \\ H^t_1 \\ H^t_2 \\ H^t_3 \\ H^t_4 
\end{bmatrix}_{(16 \times \texttt{nE}) \times 1}$ &
$\texttt{J3} = \begin{bmatrix}
    H^t_1 \\ H^t_1 \\ H^t_1 \\ H^t_1 \\ H^t_2 \\ H^t_2 \\ H^t_2 \\ H^t_2 \\ H^t_3 \\ H^t_3 \\ H^t_3 \\ H^t_3 \\ H^t_4 \\ H^t_4 \\ H^t_4 \\ H^t_4 
\end{bmatrix}_{(16 \times \texttt{nE}) \times 1}$
\end{tabular}
\]

    \item Lines \ref{line:f1}-\ref{line:f4}: We evaluate the volume forces at the centroid of faces $F\in \mathcal{F}_h$, and then we assemble
the load vector $b$ using accumarray.
\end{itemize}

\subsection{Multiplier.m}
The function \textbf{Multiplier.m} (Listing \ref{lst:LambdaPH}) implements the Lagrange multiplier matrix $C$ and vector $LN$.
\begin{lstlisting}[caption= $\mathbf{Multiplier.m}$, label=lst:LambdaPH, escapechar=|]
function [C,CT,g_ind,LN,g_NFs,Is]= Multiplier(Tetra,Coord,GI,NFs,Nb) 
nE=size(Tetra,1); 
[Fs,g_ind,~,IVRF]=fullFaces(Tetra); |\label{line:Fullfaces}|
AreaF=0.5.*vecnorm(cross(Coord(Fs(:,3),:)-Coord(Fs(:,1),:),...
    Coord(Fs(:,2),:)-Coord(Fs(:,1),:)),2,2); % area of faces
sigma=ones(4*nE,1); sigma(IVRF)=-1; |\label{line:sigma}|
r1=[1/3, 1/3 , 1/3, 0]; r2=[1/3, 0, 1/3, 1/3]; r3=[1/3, 1/3, 0, 1/3]; r4=[0, 1/3, 1/3, 1/3]; |\label{line:r1r2r3r4_Multiplier}|
CT=sigma.*AreaF.*repmat([r1;r2;r3;r4],nE,1); |\label{line:CT_Multiplier}|
I=repmat(g_ind',4,1)'; |\label{line:I_Multiplier}|
J=reshape(GI([1 5 9 13 2 6 10 14 3 7 11 15 4 8 12 16],:),4,4*nE)'; |\label{line:J_Multiplier}|
[~,g_NFs]=ismember(NFs, g_ind);  % global indices of Neumann faces  |\label{line:g_NFs}|
Is=I'; I(g_NFs,:)=[]; 
JN=J(g_NFs,:); J(g_NFs,:)=[];
CTN=CT(g_NFs,:); CT(g_NFs,:)=[];
C=sparse(I,J,CT); |\label{line:FinalC}|
%% Vectorization of Neumann BC
LN=sparse(4*nE,1); |\label{line:NeumannLN}|
if (~isempty(Nb)) |\label{line:StaterLN}|
    NP1=Coord(Nb(:,1),:); NP2=Coord(Nb(:,2),:); NP3=Coord(Nb(:,3),:);
    LNv=g(NP1,NP2,NP3).*CTN;
    LN=accumarray(JN(:),LNv(:),[4*nE,1]); 
end |\label{line:EndLN}|
end

\end{lstlisting}

\begin{itemize}
    \item Line \ref{line:Fullfaces}: The function \textbf{fullFaces.m} computes all the faces formed by \texttt{Tetra}, including repeated and unique faces. Also, the \texttt{g\_ind} is a collection of indices of all faces, with repeated face index is marked with a unique face index and size is $ (4 \times \texttt{nE}) \times 1$. The \texttt{IVRF} and \texttt{IVNRF} are index vectors of repeated and non-repeated faces, respectively.
    \item Lines \ref{line:sigma}- \ref{line:CT_Multiplier} : The basis function given  in \eqref{eq:sigma} is computed and the values of $r_1$, $r_2$, $r_3$ and $r_4$ seen in \eqref{eq:C(T)}. The \texttt{CT} stores the local Lagrange multiplier matrix $C(T)$ (see \eqref{eq:C(T)}) for each $T \in \mathcal{T}_h$.
    
    \item Lines: \ref{line:I_Multiplier}-\ref{line:J_Multiplier}: The arrays $\texttt{I}_{(4 \times \texttt{nE}) \times 4}$ and $\texttt{J}_{(4 \times \texttt{nE}) \times 4}$ are constructed form \texttt{g\_ind} and \texttt{GI} (see Listing \ref{lst:stiffPH}) simultaneously.
The entries of \texttt{I} are: \newline
$$\texttt{I} = \begin{bmatrix}\texttt{g\_ind}~~~~\texttt{g\_ind}~~~~\texttt{g\_ind}~~~~\texttt{g\_ind} \end{bmatrix}_{(4 \times \texttt{nE}) \times 4}$$

Define a row vector, $\mathcal{D}_{\beta}$=$ \begin{bmatrix}1+4\beta ~~~~ 2+4\beta ~~~~ 3+4\beta ~~~~ 4+4\beta\end{bmatrix}$, where $\beta=0,1, \cdots , \texttt{nE}-1$.

Then, using the definition of $\mathcal{D}_{\beta}$, we represent the entries of \texttt{J} as: \newline
  $$\texttt{J} = \begin{bmatrix}
   \mathcal{D}_{0} \\ \mathcal{D}_{0} \\ \mathcal{D}_{0} \\ \mathcal{D}_{0} \\ \mathcal{D}_{1} \\ \mathcal{D}_{1} \\\mathcal{D}_{1} \\ \mathcal{D}_{1} \\   \vdots \\ \mathcal{D}_{\texttt{nE}-1} \\ \mathcal{D}_{\texttt{nE}-1} \\ \mathcal{D}_{\texttt{nE}-1} \\ \mathcal{D}_{\texttt{nE}-1} 
\end{bmatrix}_{(4 \times \texttt{nE}) \times 4}.$$ 

\item \ref{line:g_NFs}-\ref{line:FinalC}: The arrays \texttt{I},~\texttt{J}, and \texttt{CT} also contain the indices and values associated with the Neumann boundary. Thus, we simultaneously compute and assemble the Lagrange multiplier matrix $C$ and the vector $LN$.

\note{The arrays $\texttt{I},~\texttt{J}$ and $\texttt{CT}$ are maintained with a size of $(4 \times \texttt{nE}) \times 4$ to allow for the easy removal or utilization of rows associated with Neumann face indices.}

 
\end{itemize}

\subsection{Schur complement and parallel computation}
The system in \eqref{eq:VectorForm} gives the Schur complement for $\Lambda$,
\begin{equation} \label{eq:schur}
    S \Lambda = -b_{D} + C \mathbb{A}^{-1} \mathbb{B}
\end{equation}
where $S_{L \times L}=-C \mathbb{A}^{-1} C',$ $ ~\mathbb{A}_{N \times N}=K+M, ~\mathbb{B}_{N \times 1}=b+LN$.
To recover $U$,
\begin{equation}
    U= \mathbb{A}^{-1} (\mathbb{B}+C'\Lambda).
\end{equation}


Algorithm \ref{alg:schur} employs a Schur complement technique to solve the linear system \eqref{eq:VectorForm} locally for each tetrahedron. By decoupling primal and hybrid variables, the method is well-suited for parallel computation in MATLAB. The \textbf{SchurComplement\_parallel.m} computes the Schur complement matrix and right-hand side locally (see lines \ref{line:localSchur}-\ref{line:localRhs}), enabling the use of \textbf{parfor} (parallel for-loop). Leveraging MATLAB’s Parallel Computing Toolbox, we distribute computations across 4 workers. Similarly, primal solution recovery also utilizes \textbf{parfor} (see lines \ref{line:StartsecondParfor}-\ref{line:EndsecondParfor}). Additionally, we vectorize the global assembly of the Schur complement matrix and right-hand side using face indices as degrees of freedom (see line \ref{line:Assembleschur}).

\algblock{ParFor}{EndParFor}
\algnewcommand\algorithmicparfor{\textbf{parfor}}
\algnewcommand\algorithmicpardo{\textbf{do}}
\algnewcommand\algorithmicendparfor{\textbf{end\ parfor}}
\algrenewtext{ParFor}[1]{\algorithmicparfor\ #1\ \algorithmicpardo}
\algrenewtext{EndParFor}{\algorithmicendparfor}
\begin{algorithm}
\footnotesize
 \caption{SchurComplement}\label{alg:schur}
\begin{algorithmic}
\State $\mathbf{Output}$: $U$, $\Lambda$
\begin{algorithmic}[1]
\ParFor{$T=1$ to \texttt{nE}} \label{line:StartFirstParfor}
\State Compute local stiffness matrix $K(T)$ \eqref{eq:K(T)}, local mass matrix $M(T)$ \eqref{eq:M(T)}, local lambda matrix $C(T)$ \eqref{eq:C(T)}, $b(T)$, $b_{D}(T)$ and $LN(T)$. 
\State Set $\mathbb{A}(T)=K(T)+M(T)$, and $\mathbb{B}(T)=b(T)+LN(T)$.
\State Compute local Schur matrix:
$$S(T)= -C(T) \mathbb{A}(T)^{-1} C(T)'$$ \label{line:localSchur} 
\State Compute local right-hand side of \eqref{eq:schur}:
$$ \texttt{rhs}(T)=-b_{D}(T) + C(T) \mathbb{A}(T)^{-1} \mathbb{B}(T) $$  \label{line:localRhs} 
\EndParFor \label{line:EndFirstParfor}
\State Assemble the global Schur matrix $S_{L \times L}$ and right hand side $\texttt{rhs}_{L \times 1}$. \label{line:Assembleschur}
\State Solve $S \Lambda = \texttt{rhs}$ \eqref{eq:schur}.
\ParFor{$T=1$ to \texttt{nE}} \label{line:StartsecondParfor}
\State Recover local primal solution:
$$U(T)=\mathbb{A}(T)^{-1} (\mathbb{B}(T)+C(T)'\Lambda(T))$$
\EndParFor   \label{line:EndsecondParfor}
\State Assemble the global primal solution $U$.
     \end{algorithmic}
\end{algorithmic}
 \end{algorithm}

\subsubsection{Vectorized implementation}
The function $\mathbf{SchurComplement\_vectorization.m}$ (Listing \ref{lst:SchurComplement}) computes the primal ($u_h$) and hybrid ($\kappa_h$) solutions.
\begin{lstlisting}[caption= $\mathbf{SchurComplement\_vectorization.m}$, label=lst:SchurComplement, escapechar=|]
function [kh,uh,size_S_global,S_global,Rhs_global]=SchurComplement_vectorization(Tetra,local_rhs,uhb1,nF,AT,CT,g_ind,g_NFs,RFaces,Is)
nE=size(Tetra,1);
CT_all=sparse(4*nE,4); |\label{line:CT_Start}|
id=setdiff(1:4*nE,g_NFs);
CT_all(id,:)=CT;   |\label{line:CT_End}|
ATall=AT(reshape(1:4*nE,nE,4)',:);  
f_all=reshape(local_rhs,[4,1,nE]); % reshape of local load vectors to 3D array
C_all=reshape(full(CT_all'),[4,4,nE]); % reshape of transpose of local lambda matrix to 3D array 
STall=-pagemtimes(pagetranspose(C_all),pagemldivide(reshape(ATall',[4,4,nE]),C_all)); |\label{line:pagemldivide}|
STall=reshape(STall, 4*nE,4);
RTall=pagemtimes(pagetranspose(C_all),pagemldivide(reshape(ATall',[4,4,nE]),f_all)); |\label{line:pagemtimes}|
RTall=reshape(RTall,4*nE,1);
%% Assemble Schur complement and rhs
Js=reshape(repmat(reshape(g_ind,4,nE),4,1),16*nE,1); |\label{line:Js}|
S_global=sparse(Is(:),Js,STall,nF,nF); |\label{line:S_global}|
S_global=S_global(RFaces,RFaces); |\label{line:S_globalRFaces}|
Rhs_global=accumarray(g_ind,RTall,[nF,1]); 
Rhs_global=Rhs_global(RFaces)-uhb1; |\label{line:Rhs_globalRFaces}|
kh=S_global \ Rhs_global; % Solve lambda
kh_all=sparse(RFaces', ones(size(RFaces')), kh, nF, 1);  %extending kappa_h
khT_all=reshape(full(kh_all(g_ind)),[4,1,nE]);  % reshape of kappa_h to 3D array
uh=pagemldivide(reshape(ATall',[4,4,nE]),f_all+pagemtimes(C_all,khT_all));% Recover uh in reverse
uh=reshape(uh,4*nE,1); 
size_S_global=size(S_global,1);
end
\end{lstlisting}
\begin{itemize}

    \item Lines \ref{line:CT_Start}-\ref{line:CT_End}: Extending the \texttt{CT} array by adding rows of zeros for the Neumann faces indices results in the $\texttt{CT\_all}$ array. 
    \item Line \ref{line:pagemldivide}: The functions \textbf{pagemtimes},  \textbf{pagetranspose}, and \textbf{pagemldivide} are used to compute local matrix multiplications, transposes, and left matrix divide, respectively. This enables vectorized operations for all $T \in \mathcal{T}_h$.
    \item Lines \ref{line:Js}-\ref{line:S_global}: We assembled the Schur complement matrix ($\texttt{S\_global}$) with respect to the degree of freedom of the face index.
    \item Lines \ref{line:S_globalRFaces}-\ref{line:Rhs_globalRFaces}: The zero rows in the $\texttt{CT\_all}$ array lead to zero rows and columns in $\texttt{S\_global}$ matrix at the indices corresponding to the Neumann faces, which were handled by the \texttt{RFaces} index vector. Similarly, the corresponding zero entries were removed from the $\texttt{Rhs\_global}$ vector.
\end{itemize}

\begin{example}
We have taken the exact solution of  \eqref{eq1} as 
\begin{equation}\label{eq:Exactsolution}
    u(x,y,z)=x^2 y^2 z^2 \qquad \text{in}~\Omega=(0,1)\times (0,1) \times (0,1)
\end{equation}
and corresponding load function $f(x,y,z)$ in the form
$$f(x,y,z)=-2(y^2z^2+x^2z^2+x^2y^2)+u(x,y,z) \qquad \text{in}~\Omega.$$
The Dirichlet and Neumann boundary conditions (cf. \figurename{~\ref{fig:3DInitialMesh}}) read

\begin{equation*}
\begin{cases}
\begin{alignedat}{2}
u_D(x,y,z) & = u(x,y,z) &\quad & \text{on } \Gamma_D,\\[6pt]
g(x,y,z) & = \bigl[\,2xy^2z^2,\; 2x^2yz^2,\; 2x^2y^2z\,\bigr]\cdot \vec{\nu}_T 
           &      & \text{on } \Gamma_N,\ \forall\, T\in\mathcal{T}_h.
\end{alignedat}
\end{cases}
\end{equation*}
The exact Lagrange multiplier from \eqref{eq1hyb1} as
 \begin{equation*}\label{eq:exactHybrid}
     \kappa=g(x,y,z)  \qquad \text{ on }~~\partial T/\Gamma_N~\forall T\in\mathcal{T}_h.
 \end{equation*}
Table~\ref{tab:StiffMultiruntime} shows the run-time (in seconds) for vectorized implementations of the functions 
$\mathbf{StiffMassPH.m}$, which computes the stiffness matrix $K_{N \times N}$, mass matrix $M_{N \times N}$, and load vector $b_{N \times 1}$, 
and $\mathbf{Multiplier.m}$, which computes the Lagrange multiplier matrix $C_{L \times N}$ and vector $LN_{N \times 1}$. 
The run-time is reported as a function of the number of tetrahedra and faces in the mesh.
The dimensions of the primal space $\mathbb{Y}_h$ (defined in \eqref{eq:Primalspace}) and the hybrid space $\mathbb{S}_h$ (defined in \eqref{eq:HybridSpace}) are
\[
N = 4 \times \texttt{nE}, \qquad 
L = n(\mathcal{F}^h_{\Omega}) + n(\mathcal{F}^h_{D}),
\]
respectively.

\begin{table}[h!]
  \centering
  \caption{Run-time of \textbf{StiffMassPH.m} and \textbf{Multiplier.m} in MATLAB.}
  \label{tab:StiffMultiruntime}
  \small\addtolength{\tabcolsep}{2pt}
  \begin{tabular}{|c|r|r|r|r|c|c|}
    \hline
   Level&  {\texttt{nE}} & \texttt{nF} & $N$ & $L$ & \multicolumn{2}{c|}{Run-time (in sec.)}   \\ \cline{6-7}
    & & & & & \textbf{StiffMassPH.m} &\textbf{Multiplier.m} \\
    \hline
  1 &  60   & 144 &240       & 104   & 0.0016 & 0.0016  \\
   2&  720 & 1,536 & 2,880     & 1376   & 0.0029  & 0.0036  \\
   3 & 8,640 & 17,664  & 34,560    & 17024  & 0.0189 & 0.0152  \\
   4&  103,680 & 208,896  & 414,720   &   206,336   & 0.3160  & 0.3358 \\
    5 &1,244,160 & 2,494,464 &  4,976,640 & 2,484,224 & 4.2434  & 4.5835  \\
    \hline
  \end{tabular}
\end{table}

Table \ref{tab:direct_vs_parallel} displays the performance comparison of the run-time (in seconds) between the direct linear solver (\textbf{mldivide}), vectorized, and a parallel Schur complement method.

\begin{table}[H]
  \centering
  \caption{Performance comparison (in seconds) between a direct solver, vectorized, and a parallel Schur complement method.}
  \label{tab:direct_vs_parallel}
  \small\addtolength{\tabcolsep}{2pt}
  \begin{tabular}{|c|r|r|r|r|r|}
    \hline
    Level & {$N$} & {$L$} & \makecell{Direct \\ solver (\textbf{mldivide})} & \makecell{Parallel Schur \\ compl. (4 workers)} & \makecell{Vectorized Schur \\ complement} \\
    \hline
    1 & 240       & 104       & 0.0011  & 0.3314  & 0.0052 \\
    2 & 2,880     & 1376     & 0.0149  & 0.1325  & 0.0083 \\
    3 & 34,560    & 17024    & 0.5833  & 0.4217  & 0.0985 \\
    4 & 414,720   & 206,336  & 17.5865 & 4.0236  & 2.5512 \\
    5 & 4,976,640 & 2,484,224 & 888.9821 & 449.5390 & 426.2138 \\
    \hline 
  \end{tabular}
\end{table}

\remark{
For the level 5, the size of the block matrix in \eqref{eq:VectorForm} is $$(N+L)\times(N+L)= 7460864 \times 7460864,$$ which results in a significantly higher runtime when solved using the default direct linear solver (888.98 seconds).  
An alternative approach is to use the Schur complement method, where the size of the Schur matrix \( S \) in \eqref{eq:schur} is 
$$L\times L=2484224 \times 2484224$$ which is approximately three times smaller than the block matrix. This reduction leads to a substantially reduced runtime of 426.21 seconds.
}

\remark{The vectorized Schur complement performs comparably to the parallel implementation when using up to four workers. With a larger number of workers, however, parallelization begins to outperform vectorization. For example, on a system equipped with an Apple M1 processor (2020) utilizing eight workers, the recorded runtimes were 206.37 seconds for the parallel Schur complement and 210.83 seconds for the vectorized Schur complement. The difference is marginal, though, since the majority of the computation time is dominated by solving the Schur linear system~\eqref{eq:schur}.
}

\medskip

In Table \ref{tab:AprioriError}, the errors $u-u_h$ in the $\mathbb{Y}$-norm, and the error $\kappa - \kappa_h$ in the $h$-norm are presented. In addition, we present the order of convergence in the respective norms.

\begin{table}[h!]
  \centering
  \caption{A priori error estimates with respect to the $\mathbb{Y}$ and $h$ -norms.}
  \label{tab:AprioriError}
  \small\addtolength{\tabcolsep}{2pt}
  \begin{tabular}{|c|c|c|c|c|}
    \hline
    Level & $\|u - u_h\|_{\mathbb{Y}}$ & O.C. ($\mathbb{Y}$) & $\|\kappa - \kappa_h\|_{h}$ & O.C. ($h$)\\
    \hline
    1         & 0.0606  & -- & 0.1790      & --\\    
   2    &   0.0303  & 1.0013  &    0.0938      &  0.9315         \\
   3  &   0.0154    &   0.9741  & 0.0487     &  0.9452    \\
   4   &  0.0078  &   0.9754 & 0.0252      &  0.9501      \\
    5 &    0.0040 & 0.9846 &     0.0131 &  0.9439   \\
    \hline
  \end{tabular}
\end{table}

The \figurename{~\ref{level4PrimalSolution}} visualizes the exact solution \eqref{eq:Exactsolution} along with the approximate solution $u_h$ at mesh refinement level 4. Additionally, \figurename{~\ref{level4ExtHybridSolution}} \& \ref{level4IntHybridSolution} depict the exact and approximate Lagrange multiplier on the faces belonging to $\mathcal{F}^h_{D}$ and $\mathcal{F}^h_{\Omega}$.

\begin{figure}
\vspace{-2cm}
\centering 
\begin{subfigure}
 \centering 
 \includegraphics[scale=0.34]{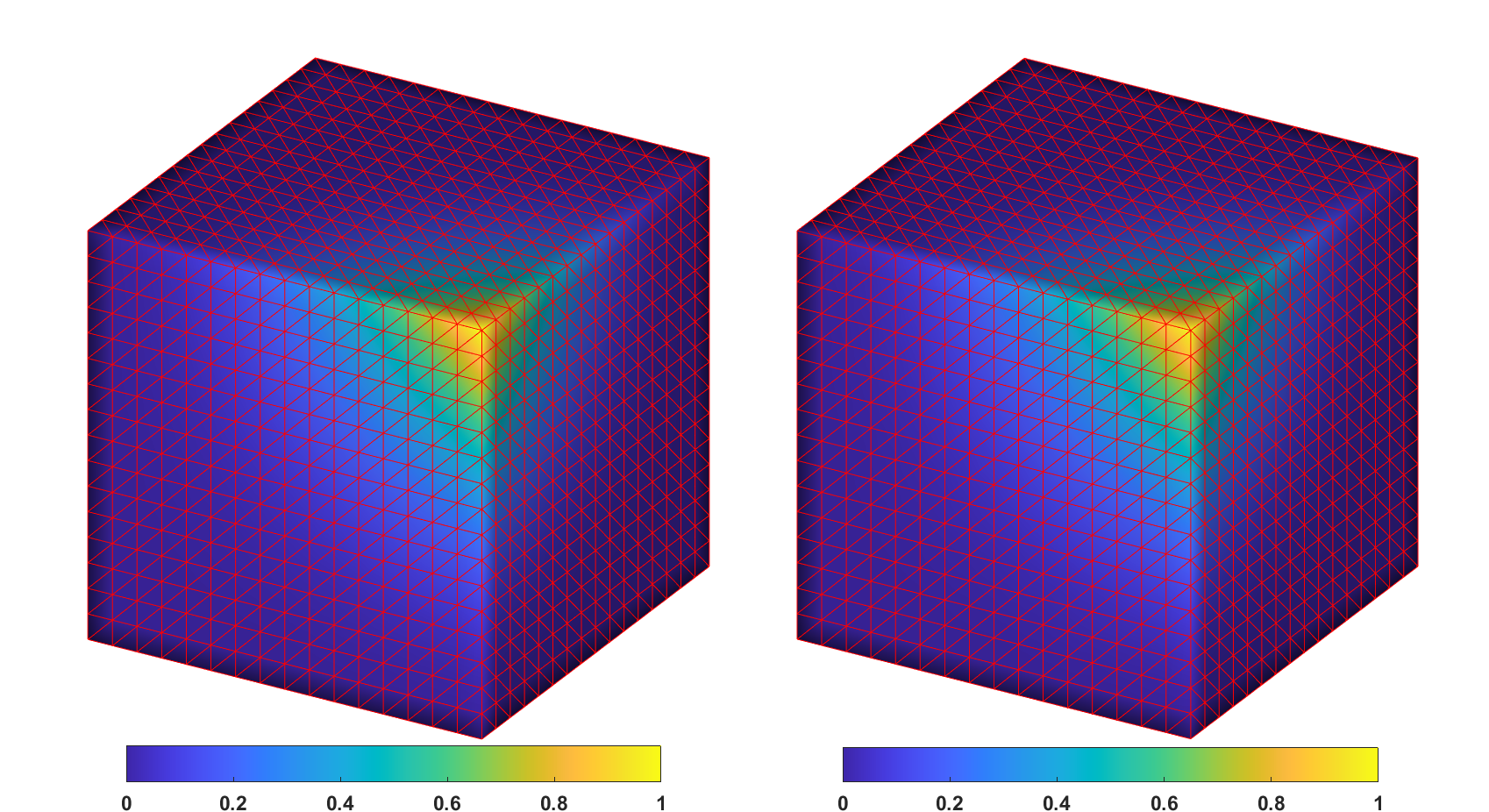}
 \vspace{-0.1cm}
 \caption{Visualization of $u$ (left) and  $u_h$ (right)  at Level 4.}
 \label{level4PrimalSolution}
\end{subfigure}
\vspace{0.4cm}
\begin{subfigure}
 \centering 
 \includegraphics[scale=0.48]{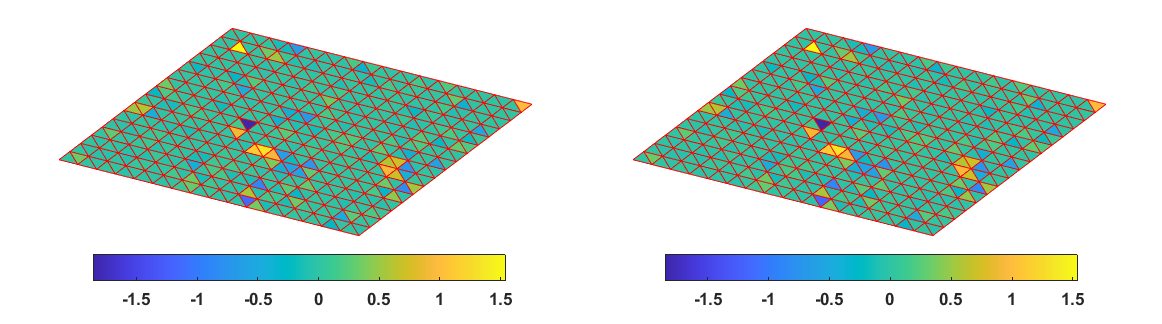}
  \vspace{-0.4cm}
 \caption{Visualization of $\kappa$ (left) and  $\kappa_h$ (right)  on exterior faces in $\mathcal{F}^h_{D}$ at Level 4.}
 \label{level4ExtHybridSolution}
\end{subfigure}
\begin{subfigure}
 \centering 
 \includegraphics[scale=0.345]{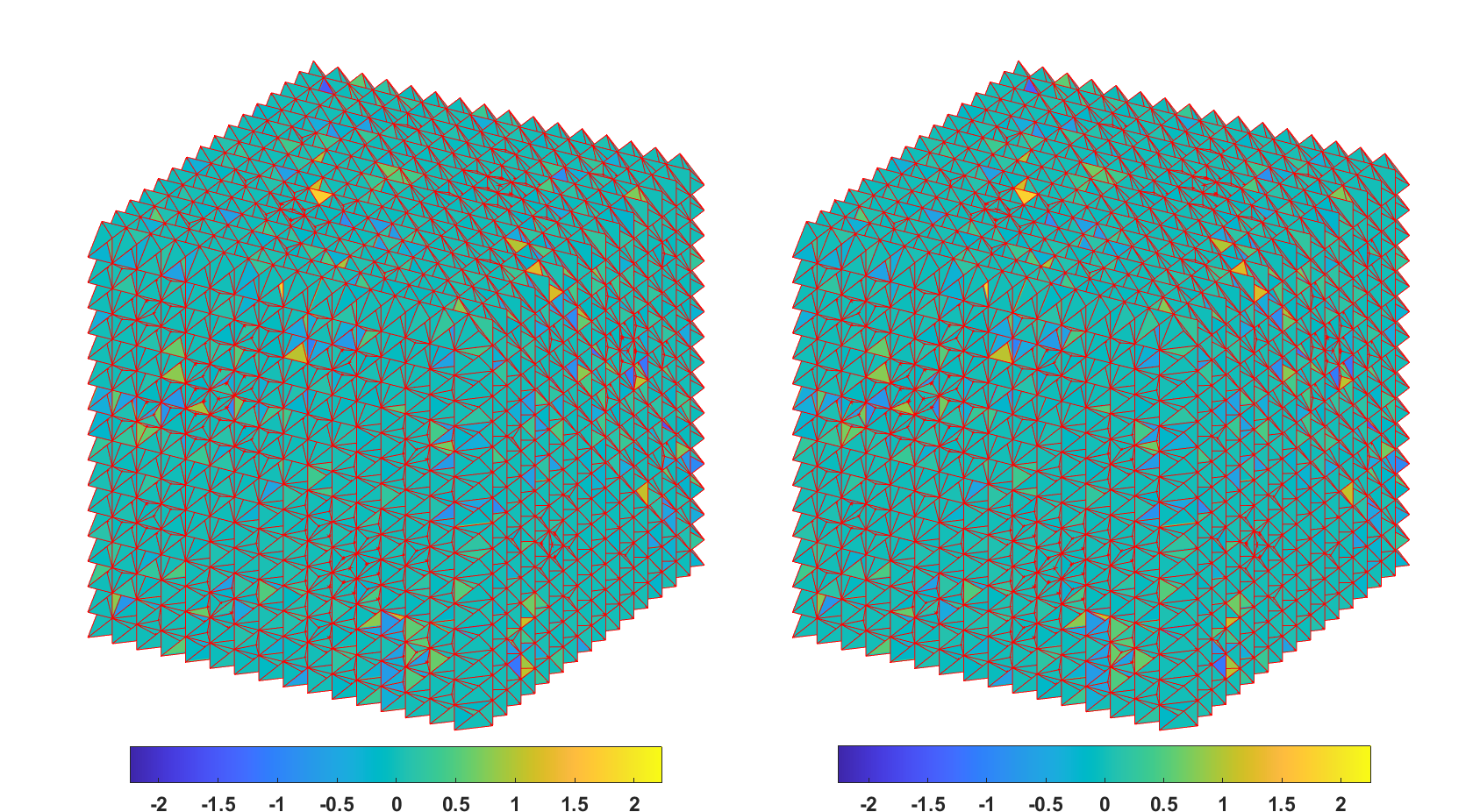}
  \vspace{-0.2cm}
\caption{Visualization of $\kappa$ (left)  and  $\kappa_h$ (right) on interior faces in $\mathcal{F}^h_{\Omega}$ at Level 4.}
  \label{level4IntHybridSolution}
\end{subfigure}
\end{figure}

\end{example}

\section{Conclusion}\label{sub:conclusion}
The 3D uniform mesh refinement and primal hybrid FEM for elliptic problems with mixed boundary conditions have been efficiently implemented in MATLAB. The simultaneous computation of primal and
hybrid solutions without the requirement of inter-element boundaries allows successful parallel computation using the Schur complement method. The numerical results show the performance of our MATLAB software package and the optimal order of convergence of the method. All codes are flexible and easily adaptable to other programming languages. In future work, we will include the $hp$-adaptivity of the primal hybrid FEM, also for quadrilateral-type elements. 

\section*{Acknowledgement}

  The third author announces the support of the Czech Science Foundation (GACR) through the GA23-04766S grant, Variational approaches to dynamical problems in continuum mechanics.
\section*{Declarations}

\noindent
\textbf{Conflict of interests}: 
The authors have no known conflicts of financial interests or personal relationships that could have appeared to influence the work reported in this manuscript.

 \bibliography{Main}

 \newpage

\end{document}